\documentclass[11pt]{article}
\usepackage[english]{style}
\usepackage{latexsym,amsmath,amssymb,amsfonts,diagrams,graphics,amscd}
\usepackage{graphicx}
\input xypic
\xyoption{all}
\input epsf.tex

\usepackage{pgf,tikz}
\usetikzlibrary{arrows}
\DeclareFontFamily{U}{rsf}{}
\DeclareFontShape{U}{rsf}{m}{n}{
  <5> <6> rsfs5 <7> <8> <9> rsfs7 <10-> rsfs10}{}
\DeclareMathAlphabet{\mathscr}{U}{rsf}{m}{n}

\DeclareMathAlphabet{\mathgth}{U}{euf}{m}{n}

\DeclareFontFamily{U}{cyr}{}
\DeclareFontShape{U}{cyr}{m}{n}{
  <5> wncyr5 <6> wncyr6 <7> wncyr7 <8> wncyr8 <9> wncyr9 <10-> wncyr10}{}
\DeclareMathAlphabet{\mathcyr}{U}{cyr}{m}{n}

\input cyracc.def

\makeatletter
\def\operator@font{\sf}
\makeatother

\newcommand{\sbt}{{\scalebox{0.5}{\textbullet}}}
\newcommand{\cA}{{\mathscr A}}
\newcommand{\cB}{{\mathscr B}}
\newcommand{\sA}{{\mathcal A}}

\newcommand{\sL}{{\mathcal L}}

\newcommand{\cG}{{\mathscr G}}
\newcommand{\coH}{{\mathscr H}}
\newcommand{\cI}{{\mathscr I}}

\newcommand{\cO}{{\mathscr O}}
\newcommand{\OO}{{\mathcal O}}

\newcommand{\bN}{{\overline{N}}}

\newcommand{\bX}{{\overline{X}}}
\newcommand{\bY}{{\overline{Y}}}
\newcommand{\bW}{{\overline{W}}}
\newcommand{\bS}{{\overline{S}}}

\newcommand{\bi}{{\bar{i}}}
\newcommand{\bj}{{\bar{j}}}

\newcommand{\tF}{\tilde{F}}
\newcommand{\tf}{\tilde{f}}

\newcommand{\tX}{\tilde{X}}

\newcommand{\dR}{{\mathsf dR}}

\newcommand{\sHom}{\underline{\mathsf{Hom}}}
\newcommand{\RsHom}{{\R\sHom}}

\newcommand{\sEnd}{\underline{\mathsf{End}}}

\newcommand{\D}{{\mathbf D}}

\newcommand{\chk}{{\scriptscriptstyle\vee}}
\newcommand{\R}{\mathbf{R}}

\newcommand{\bbk}{\mathbf{k}}

\newcommand{\bbE}{\mathbb{E}}

\newcommand{\Gm}{{\mathbb G}_m}
\renewcommand{\S}{{\mathbb{S}}}

\DeclareMathOperator{\Tors}{Tors}
\DeclareMathOperator{\Split}{Spl}
\DeclareMathOperator{\Spec}{Spec}

\DeclareMathOperator{\End}{End}

\DeclareMathOperator{\Tot}{Tot}

\DeclareMathOperator{\Pic}{Pic}

\DeclareMathOperator{\id}{id}

\DeclareMathOperator{\rk}{rk}
\DeclareMathOperator{\Ext}{Ext}

\DeclareMathOperator{\codim}{codim}

\DeclareMathOperator{\Crit}{Crit}
\DeclareMathOperator{\AH}{hexp}

\newcommand{\ra}{\rightarrow}

\newcommand{\lra}{\longrightarrow}

\newcommand{\inj}{\hookrightarrow}

\newcommand{\C}{\mathbb{C}}

\newcommand{\gMod}{\mathgth{Mod}}
\newcommand{\gCoh}{\mathgth{Coh}}

\newcommand{\iso}{\cong}

\newcommand{\del}{\partial}

\newcommand{\EG}{\mathcal{EG}}
\newcommand{\field}[1]{\mathbb{#1}}

\newcommand{\A}{\field{A}}
\renewcommand{\phi}{\varphi}

\newarrow{Equal} =====

\author{%
Dima Arinkin\thanks{Mathematics Department,
University of Wisconsin--Madison, 480 Lincoln Drive, Madison, WI
53706, USA, {\em e-mail: }{\tt arinkin@math.wisc.edu,
andreic@math.wisc.edu}},\ \ 
Andrei C\u ald\u araru,\footnotemark[1]\ \  M\'arton
Hablicsek\thanks{Mathematics Department, University of Pennsylvania,
David Rittenhouse Lab, 209 S.\ 33rd Street, Philadelphia, PA 19104,
USA, {\em e-mail: }{\tt mhabli@math.upenn.edu}}
}

\title{Derived intersections and the Hodge theorem}

\begin{document}

\maketitle

\begin{abstract} 
  The algebraic Hodge theorem was proved  in a beautiful 1987 paper by Deligne and Illusie, using positive characteristic methods. We argue that the central algebraic object of their proof can be understood geometrically as a line bundle on a derived scheme.  In this interpretation, the Deligne-Illusie result can be seen as a proof that this line bundle is trivial under certain assumptions.

  We give a criterion for the triviality of this line bundle in a more general context.  The proof uses techniques from derived algebraic geometry, specifically arguments which show the formality of certain derived intersections.  Applying our criterion we recover Deligne and Illusie's original result.  We also apply these techniques to the result of Barannikov-Kontsevich, Sabbah, and Ogus-Vologodsky concerning the formality of the twisted de Rham complex.
\end{abstract}

\section*{Introduction}
In their 1987 paper~\cite{DelIll}, Deligne and Illusie proved the
following algebraic version of the Hodge theorem.

\begin{Theorem}
\label{DelIll}
Let $X$ be a smooth proper scheme over a perfect field $\bbk$ of positive
characteristic $p>\dim X$.  Assume that $X$ lifts to the ring $W_2(\bbk)$
of second Witt vectors of $k$.  Then the Hodge-to-de Rham spectral
sequence for $X$ degenerates at ${}^1 E$.
\end{Theorem}

\medskip

\noindent
The corresponding result in characteristic zero (the Hodge theorem)
follows easily from this result using a standard argument of
Frobenius.

\paragraph
\label{mainresult}
Deligne and Illusie deduce Theorem~\ref{DelIll} from the following
statement, whose proof forms the central part of their paper.  Consider the
Frobenius twist $X'$ of $X$, and the corresponding relative Frobenius
morphism $F:X\ra X'$.  If we denote by $\Omega_X^\sbt$ the algebraic
de Rham complex of $X$, then the complex $F_*\Omega_X^\sbt$ is an
object of $\D(X')$, the derived category of coherent sheaves on $X'$.
Deligne and Illusie proved that the two statements below are
equivalent.
\begin{itemize}
\item[(1)] $X$ lifts to $W_2(\bbk)$.
\item[(2)] $F_*\Omega^{\sbt}_X$ is formal in $\D(X')$.
\end{itemize}
(A complex is said to be {\em formal} if it is isomorphic to the direct sum of
its cohomology sheaves.)
\medskip

\noindent
Independently of the work of Deligne and Illusie the first two authors
proved the following result in~\cite{AriCal}.

\begin{Theorem}
\label{AriCal}
Let $i:X\inj S$ be a closed embedding of smooth schemes over a field
$\bbk$ of characteristic $p>\codim_S X$, and assume that the normal
bundle $N = N_{X/S}$ is the restriction of a vector bundle $\bN$ on
the first infinitesimal neighborhood of $X$ in $S$. Then the two
statements below are equivalent.
\begin{itemize}
\item[(1)] The trivial vector bundle $\cO_X$ lifts to the first
  infinitesimal neighborhood of $X$ in $S$.
\item[(2)] The structure complex $i^*i_*\cO_X$ of the derived
  self-intersection of $X$ inside $S$ is formal.  
\end{itemize}
(All functors are implicitly assumed to be derived.)
\end{Theorem}

\paragraph {\bf Remark.} As stated, condition (1) is vacuous.  We'll
see that in order to understand the result of Deligne and Illusie we
need to work with a generalization of Theorem~\ref{AriCal} to Azumaya
spaces (see below), where (1) becomes a real condition.

\paragraph
Mircea Musta\c t\u a noted that there is a strong similarity between
Theorem~\ref{AriCal} and the Deligne-Illusie equivalence.  He asked if
it is possible to recast the problem of Deligne-Illusie as a derived
intersection problem such that their result would follow from
Theorem~\ref{AriCal}.  The goal of this paper is to give a positive
answer to this question, by constructing a closed embedding $h:X' \inj
\bS$ for which $F_*\Omega_X^\sbt$ is dual to the structure sheaf
$h^*h_* \cO_{X'}$ of the derived self-intersection of $X'$ inside
$\bS$.  In particular $F_*\Omega_X^\sbt$ will be formal if and only if
$h^*h_*\cO_{X'}$ is.

One feature of this construction is that $\bS$ is not an ordinary
scheme, but an extended object called an {\em Azumaya space}, see
Section~\ref{sec:Azumaya} or~\cite{CalThesis}.  Such a space is a pair
$(S, \cA)$ where $S$ is an ordinary scheme and $\cA$ is a sheaf
of Azumaya algebras on $S$.  Coherent sheaves on the Azumaya
space $(S, \cA)$ are coherent sheaves on $S$ together with a left
$\cA$-module structure.  Morphisms in the category of Azumaya spaces
(or, more precisely, 1-morphisms in this 2-category) are given by
Morita equivalences.

Because of how morphisms are defined in the category of Azumaya
spaces, it is not {\em a priori} true that the the trivial bundle
$\cO_{X'}$ extends to the first infinitesimal neighborhood for the
embedding $h$.  In fact we will prove that in the Deligne-Illusie
problem this happens if and only if $X$ lifts to $W_2(\bbk)$.  Applying a
version of Theorem~\ref{AriCal} suited to Azumaya spaces
(Corollary~\ref{twisted splitting}) in the setting we have engineered
above recovers Theorem~\ref{DelIll}.

\paragraph
We now describe the construction of the embedding $h$.  For this we
need some classic facts about rings of crystalline differential
operators (for more details see Section~\ref{sec:basics}).  Denote the
total space of the cotangent bundle of $X'$ by $T^*X'$, and consider
the embedding $i:X'\inj T^*X'$ of $X'$ into $T^*X'$ as the zero
section. The sheaf of crystalline differential operators $D_X$ of $X$
has center $\cO_{T^*X'}$, and thus can be regarded as a sheaf of
algebras $D$ over $T^*X'$. Moreover $D$ is an Azumaya algebra over
$T^*X'$, and along the zero section of $T^*X'$ the Azumaya algebra $D$
splits,
\[ D|_{X'} \iso \End_{X'}(F_*\cO_X). \] 
Thus $\cO_{X'}$ and $D|_{X'}$ are Morita equivalent, and we obtain an equivalence
of abelian categories
\[ \gCoh(X') \iso \gCoh(X', D|_{X'}), \] 
where the latter is the abelian category of coherent sheaves on $X'$
endowed with the structure of $D|_{X'}$ left module.  

\paragraph
From the point of view of Azumaya spaces the above equivalence of
categories of coherent sheaves is simply an isomorphism
\[ m: (X',\cO_{X'}) \xrightarrow{\sim} (X', D|_{X'}). \] 
Moreover, the embedding $i:X'\ra T^*X'$ gives rise to an embedding of
Azumaya spaces
\[ i_D : (X', D|_{X'}) \ra (T^*X', D). \]
The first result of this paper is the following.

\begin{Theorem}
\label{deRham}
Let $X$ be a smooth scheme over a perfect field $\bbk$ of characteristic $p>0$.
Consider the composite morphism
\[ h:(X',\OO_{X'})\xrightarrow{m}(X',D|_{X'})\xrightarrow{i_D}(T^*X',D). \]
Then there are natural isomorphisms in $\D(X')$
\[ F_*\Omega^{\sbt}_X\iso h^!h_*\OO_{X'} \iso \left ( h^*h_* \cO_{X'}
\right )^\chk. \]
Thus $F_*\Omega^\sbt_X$ is formal if and only if $h^*h_*\cO_{X'}$ is.
\end{Theorem}

\paragraph {\bf Twisted derived intersections.}
\label{subsec:twiderint}
This theorem shows that it is important to understand formality
properties for derived intersections inside Azumaya spaces.  This is a
more general problem which we discuss now.

Let $\bS = (S, \cA)$ be an Azumaya space, and let $X$ and $Y$
be subvarieties of $S$ such that $\cA$ splits along $X$ and $Y$.  By
choosing splitting modules, one can regard the ordinary (untwisted)
spaces $X$, $Y$ as subvarieties of $(S, \cA)$.  We want to compare the
dg spaces $W$, $\bW$ which are obtained as the derived intersections
of $X$ and $Y$ inside $S$ and $\bS$, respectively.

\paragraph
The central geometric observation of this paper is that there is a
natural line bundle $\sL$ on $W$ which measures the difference between
the spaces $W$ and $\bW$.  (In a sense, this line bundle measures
the extent to which the chosen splittings fail to agree along $W$, see
Section~\ref{sec:twiderint}.)  If the dg space $W$ is formal,
formality of $\bW$ often reduces to the problem of checking whether $\sL$ is
trivial.

\paragraph
This discussion allows us to recast the original Deligne-Illusie
problem in geometric language.  For this problem we need to study the
derived self-intersection of the zero section $X'$ inside the Azumaya
space $(T^*X', D)$.  The untwisted intersection $W$ is formal (as the
self-intersection of the zero section in a vector bundle).
Theorem~\ref{deRham} shows that the line bundle $\sL$ on $W$
associated to this twisted intersection problem is
$(F_*\Omega^\sbt_X)^\chk$.  This gives a geometric meaning to the
obstruction $\alpha$ to the formality of $F_*\Omega^\sbt_X$: it is
precisely the class of $\sL$ in $\Pic(W)$.  Indeed, in
Deligne-Illusie's paper it is shown that $\alpha$ lies in $H^2(X',
T_{X'})$, and we'll argue that this group is a direct summand in
$\Pic(W)$, thus confirming the geometric description above.

\paragraph
One can ask for a criterion, for an arbitrary twisted derived
intersection problem, that will imply the triviality of the associated
line bundle $\sL$.  We give such a criterion in
Theorem~\ref{thm:twiderint}.  For a self-intersection this result is
similar in nature to Theorem~\ref{AriCal}; for an arbitrary
intersection it is a generalization to Azumaya spaces of the main
result of~\cite{AriCalHab}.

For the Deligne-Illusie application, the main condition of
Theorem~\ref{thm:twiderint} that needs to be satisfied in order to
conclude the formality of $F_*\Omega_X^\sbt$ is for the trivial bundle
of $X'$ to extend to the first infinitesimal neighborhood of the
embedding $h$.  This extension problem needs to be understood
correctly in the context of maps of Azumaya spaces.  We'll argue (see
proof of Theorem~\ref{twisted}) that the trivial bundle extends if and
only if the Azumaya algebra $D$ splits on the first infinitesimal
neighborhood of the zero section in $T^*X'$.  \medskip

\noindent
In this context we can state the second result of this paper.
\medskip

\begin{Theorem}
\label{statements}
Let $X$ be a smooth scheme over a perfect field $\bbk$ of characteristic
$p>\dim X$.  Then the following five statements are equivalent. 
\begin{itemize}
\item[(1)] $X$ lifts to $W_2(\bbk)$.
\item[(2)] $D$ splits on the first infinitesimal neighborhood of $X'$ in $T^*X'$.
\item[(3)] the trivial bundle $\cO_{X'}$ lifts to the
  first infinitesimal neighborhood of the embedding $h$.
\item[(4)] the associated line bundle $\sL\in \Pic(W)$ is trivial.
\item[(5)] $F_*\Omega^{\sbt}_X$ is formal in $\D(X')$.
\end{itemize}
\end{Theorem}
\medskip

\paragraph
{\bf Remark.}
Some of the implications above are not new.  The equivalence of $(1)$
and $(5)$ is the original result of Deligne and
Illusie~\cite{DelIll}. The implication $(1)\Rightarrow (2)$ was proven
by Ogus and Vologodsky~\cite{OguVol}, where a splitting module is
explicitly constructed. The implications $(2) \Rightarrow (3)
\Rightarrow (4) \Rightarrow (5) \Rightarrow (2)$ are parts of
Theorem~\ref{thm:twiderint}, and give a geometric interpretation of
Deligne and Illusie's result.

\paragraph
{\bf The twisted de Rham complex.}
There is a generalization of the de Rham theorem that is of interest
in the study of singularity theory and in the theory of matrix
factorizations, and which involves {\em twisted de Rham complexes}.
This generalization was investigated by Barannikov and
Kontsevich~\cite{BarKon} and Sabbah~\cite{Sab}, and in the remainder
of this introduction we discuss how some of their results can also be
understood as formality results for derived intersections.

\paragraph 
Let $X$ be a smooth quasi-projective scheme over $\C$, with a map
$f:X\ra \A^1$ such that its critical locus is proper over $\A^1$. These
data give rise to two complexes defined as
\begin{align*}
\Omega^{\sbt}_{X,d-\wedge df} & : \quad 0\rightarrow \OO_X\xrightarrow{d-\wedge df}\Omega^1_X\xrightarrow{d-\wedge
df}\dots \xrightarrow{d-\wedge df}\Omega^n_X\rightarrow 0,
\intertext{and}
\Omega^{\sbt}_{X,\wedge df} & : \quad 0\rightarrow \OO_X\xrightarrow{\wedge df}\Omega^1_X\xrightarrow{\wedge
df}\dots \xrightarrow{\wedge df}\Omega^n_X\rightarrow 0.
\end{align*}

\noindent
The following theorem, due to Barannikov-Kontsevich~\cite{BarKon}
and Sabbah~\cite{Sab}, relates the hypercohomology of these two complexes.

\begin{Theorem}\label{BarKon}
The hypercohomology spaces of the two complexes defined above are of the same
(finite) dimensions.
\end{Theorem}

\paragraph
Note that the de Rham theorem is the particular case of the above
result when $f = 0$.  In~\cite{Sab} Sabbah asked for a purely
algebraic proof of the result of Barannikov-Kontsevich.  Such a proof
 was given by Ogus and
Vologodsky~\cite{OguVol}, again using positive characteristic methods.

We can ask again if these results can be rephrased as
formality statements for derived intersections, and in what follows we
shall give a generalization of Theorem~\ref{deRham} and a partial
version of the equivalence in Theorem~\ref{statements}. 

\paragraph
Consider a smooth scheme $X$ over a perfect field $\bbk$ of positive
characteristic $p$, with a regular function $f:X\ra \A^1$. The map
$f$ provides a morphism $f':X'\ra \A^1$.  We denote the graph of
$d(f')$ in $T^*X'$ by $X'_f$, and it is easy to see~(\ref{def:L})
that $D$ splits on the section $X'_f$ just as it did on the zero
section $X'$, giving rise to an isomorphism
\[ n:(X'_f, \cO_{X_f'}) \ra (X'_f, D|_{X_f'}).\]
The following result generalizes Theorem~\ref{deRham}.

\begin{Theorem}
\label{twisteddeRham}
Consider the embeddings $i:X'\inj T^*X'$ and $j:X'_f\inj
T^*X'$, and the composite morphisms
\begin{align*}
i' & :(X',\OO_{X'})\xrightarrow{m}(X',D|_{X'})\xrightarrow{i_D}(T^*X',D)
\intertext{and}
j' &
:(X'_f,\OO_{X'_f})\xrightarrow{n}(X'_f,D|_{X'_f})\xrightarrow{j_D}(T^*X',D).
\end{align*}
Then we have isomorphisms in $\D(X')$
\begin{align*}
\Omega^{\sbt}_{X',\wedge df'} & \iso i^!j_*\OO_{X'_f}, \\
F_*\Omega^{\sbt}_{X,d-\wedge df} &\iso i'{}^!j'_*\OO_{X'_f}.
\end{align*}
\end{Theorem}

\paragraph
Note that we are precisely in the general context of twisted derived
intersections set up in~(\ref{subsec:twiderint}): we are studying the
derived intersections $W$, $\bW$ of $X'$ and $X'_f$ inside $S = T^*X'$
and $\bS = (T^*X', D)$, respectively.  Indeed, the above theorem shows
that $\cO_W$ and $\cO_\bW$ are precisely the duals of
$\Omega^{\sbt}_{X',\wedge df'} $ and $F_*\Omega^{\sbt}_{X,d-\wedge
  df}$.
\medskip

\noindent
Now we can apply the full power of Theorem~\ref{thm:twiderint} to
prove the following partial version of the Barannikov-Kontsevich
claim.

\begin{Theorem}
\label{thm:twisteddR}
Let $X$ be a smooth scheme over a field of characteristic zero.  Let
$f$ be a regular function on $X$ so that the map $f:X\ra \A^1$ is
proper, and assume that the following two conditions hold.
\begin{itemize}
\item[(1)] The critical locus $Z = \Crit f$ is scheme-theoretically smooth.
\item[(2)] The embedding $p:Z \hookrightarrow X$ is split to first order, that
is, the short exact sequence
\[ 0 \ra T_Z \ra T_X|_Z \ra N_{Z/X} \ra 0 \]
splits.
\end{itemize}
Then the complexes $F_*\Omega_{X, d-\wedge
df}^\sbt$ and $\Omega_{X', \wedge df'}^{\sbt}$ are both formal and there
exist isomorphisms
\[ \R^i\Gamma(X, \Omega_{X, d-\wedge df}^\sbt) \iso \R^i\Gamma(X,
\Omega^{\sbt}_{X, \wedge df}) \iso \bigoplus_{i-c=p+q} H^p(Z,\Omega^q_Z\otimes \omega),\] 
where $c$ denotes the codimension of $Z$ in $X$ and $\omega$ denotes
the relative dualizing complex of the embedding $Z\inj X$. 
\end{Theorem}  

\paragraph 
The paper is organized as follows. In Section~\ref{sec:basics} we
collect basic facts about schemes over fields of positive
characteristic and about the sheaf of crystalline differential
operators.  We follow with Section~\ref{sec:Azumaya} where we set up
the 2-category of Azumaya spaces.  In Section~\ref{sec:twiderint} we
discuss the general problem of twisted derived intersections and we
prove Theorem~\ref{thm:twiderint}.  We conclude the paper with a final
section where we apply the general results we have developed to the
study of the formality problem for de Rham complexes.  In particular we
prove Theorems~\ref{deRham},~\ref{twisteddeRham} and~\ref{thm:twisteddR}.

\paragraph 
\textbf{Acknowledgements.}  The entire project described in this paper
was sparked by a comment/question by Mircea Musta\c t\u a, for which
we are grateful.  The second author has benefited from a stimulating
conversation with Michael Groechenig.  The idea of Remark~\ref{subsec:defh} 
arose from
conversations with Taylor Dupuy.  The authors are supported by
the National Science Foundation under Grants No.\ DMS-0901224,
DMS-1101558, and DMS-1200721.

\section{Background on schemes over fields of positive characteristics}\label{sec:basics}

\paragraph
\label{frobenius}
Let $X$ be a smooth scheme over a perfect field $\bbk$ of characteristic
$p>0$. The classical Frobenius morphism on $X$ acts trivially on the
underlying topological space, but acts by the $p$-th power map on the
structure sheaf. The Frobenius twist $X'$ of $X$ is obtained by taking
the base change of $X$ by the classical Frobenius map $\Spec
\bbk\rightarrow \Spec \bbk$. Since the classical Frobenius morphisms on $X$
and on $\Spec \bbk$ commute with the structure map $X\rightarrow \Spec \bbk$
we obtain a commutative diagram
$$\xymatrix{X \ar[rd] \ar[r]^{F} & X' \ar[d] \ar[r]^{\pi} & X \ar[d]\\
&\Spec \bbk \ar[r] & \Spec \bbk.}$$
The map $F:X\ra X'$ in the above diagram is called the {\em relative
  Frobenius} morphism.  By construction it is a map of schemes over
$\bbk$.  

Since $\bbk$ is a perfect field, the map $\Spec \bbk\ra \Spec \bbk$ is an isomorphism.
Therefore, the map $\pi:X'\ra X$ is an isomorphism
of abstract schemes. However, this isomorphism is \emph{not} a morphism
of schemes over $\bbk$.

\paragraph
The algebraic de Rham complex of $X$ is defined as the complex
\[ \Omega_{X}^{\sbt} = 0 \ra \Omega^0_X \xrightarrow{d} \Omega^1_X
\xrightarrow{d}\cdots\rightarrow 0,\] 
where $d$ is the usual de Rham differential on algebraic
forms.   The algebraic de Rham cohomology of $X$ is defined as
\[ H^*_\dR(X) = R^*\Gamma(X, \Omega_X^\sbt). \] 
The Hodge-de Rham spectral sequence 
\[ {}^1 E^{pq} = H^p(X, \Omega^q_X) \Rightarrow H^{p+q}_\dR(X)\]
arises from the ``stupid'' filtration of $\Omega_X^\sbt$, whose
associated graded terms are $\Omega^q_X$.  

It is important to note that while the sheaves in the complex
$\Omega_X^\sbt$ are coherent $\cO_X$-modules, the differentials are
{\em not} $\cO_X$-linear, so we shall only consider the above complex
as a complex of sheaves of abelian groups. On the other hand the equality
\[ d(f^pg)=d(f^p)\cdot g+f^p\cdot d(g)=pf^{p-1}df\cdot g+ f^p\cdot
d(g) \]
implies that $F_*\Omega^{\sbt}_X$ is a complex of coherent
$\OO_{X'}$-modules (the differentials are now $\OO_{X'}$-linear).

\paragraph
The main result of Deligne-Illusie, Theorem~\ref{DelIll}, asserts that
if $X$ is proper over $\bbk$ and lifts to $W_2(\bbk)$ then the Hodge-to-de Rham
spectral sequence degenerates at ${}^1 E$.  This follows easily from
the central statement that under the liftability assumption to
$W_2(\bbk)$ we have $F_*\Omega_X^\sbt$ is formal.  We describe this
implication.

A well-known result of Cartier describes the cohomology sheaves of
$F_*\Omega_X^\sbt$ for any $X$ smooth over $k$:
\[ \mathcal{H}^i(X',F_*\Omega^{\sbt}_X)\iso \Omega_{X'}^i. \]
If $X$ lifts to $W_2(\bbk)$ then $F_*\Omega_X^\sbt$ is formal, therefore
\[F_*\Omega^{\sbt}_X\iso \bigoplus_{i=0}^{\dim X} \Omega^i_{X'}[-i].\] 
Taking hypercohomology of both sides we obtain an isomorphism
\[H_{\dR}^i(X)\iso\bigoplus_{i=p+q} H^p(X',\Omega^q_{X'}) \iso
\bigoplus_{i=p+q} H^p(X,\Omega^q_{X}),\] 
where the latter isomorphism is a consequence of the fact that $X$ and
$X'$ are isomorphic as abstract schemes.  The assumption that $X$ is
proper over $k$ implies that the vector spaces $H^p(X, \Omega_X^q)$
are finite dimensional, and the degeneration at ${}^1 E$ of the
Hodge-to-de Rham spectral sequence follows by comparing dimensions.

\paragraph
Another key ingredient in what follows is the sheaf of crystalline
differential operators (see~\cite{BezMirRum} for a detailed
description).  We emphasize that we do {\em not} work with
Grothendieck's ring of differential operators $\mathbb{D}$.  Instead
we use the ring of crystalline differential operators $D_X$, which is
the quasi-coherent sheaf of algebras over $X$ defined as the universal
enveloping algebra of the Lie algebroid of the tangent bundle $T_X$.
Explicitly, it is defined locally as the $k$ algebra generated by
$T_X$ and $\OO_X$ subject to the relations
\begin{itemize}
\item $f\cdot \partial=f\partial$ for all $f\in \OO_X$ and $\partial\in T_X$,
\item $\partial_1\cdot \partial_2-\partial_2\cdot\partial_1=[\partial_1,\partial_2]$ for all $\partial_i\in T_X$, and
\item $\partial\cdot f-f\cdot\partial=\partial(f)$ for all $f\in \OO_X$ and $\partial\in T_X$.
\end{itemize}
In the case of $X=\A^n = \Spec k[\underline{x}_i]$, this construction
produces the Weyl algebra $D_X=k\langle
\underline{x}_i,\underline{\partial}_i\rangle$. 

The main difference between $D_X$ and $\mathbb{D}$ is best understood
in the case where $X = \A^1 = \Spec k[x]$.  If we let $\del =
\del/\del x$, then in the ring $\mathbb{D}$
the operator $\del^p$ is zero because its action on 
functions is trivial.  In the ring $D_X$, $\del^p$ is non-trivial.  On
the other hand, in $\mathbb{D}$ there is an element of the form
$\del^p/p!$ which does not exist in $D_X$. 

\paragraph
Explicitly, a $D_X$-module $F$ is an $\cO_X$-module together with an action of the Lie algebroid $T_X$. Such an action can be 
represented as an integrable connection
\[\nabla=\nabla_F:F\ra\Omega_X^1\otimes F.\]
In local coordinates $\{x_1,...,x_n\}$, the connection is given by
\[\nabla(f)=\sum_{i=1}^n dx_i\otimes \frac{\del}{\del{x_i}}(f).\] 

\paragraph\label{deRham as Ext}
The sheaf $\cO_X$ is naturally a (left) $D_X$-module. This $D_X$-module admits a Koszul resolution
of the form
\[0\to D_X\otimes(\wedge^{\dim X}T_X)\to\dots\to D_X\otimes T_X\to
D_X\to\cO_X \to 0.\] 
Here the tensor products are over $\cO_X$, with $\cO_X$ acting on $D_X$ on the right.
This is a resolution of $\cO_X$ by locally free $D_X$-modules; therefore, it can be used 
to compute local Ext's. Accordingly, if $F$ is any $D_X$-module, the derived sheaf of homomorphisms
$\RsHom_{(X, D_X)}(\cO_X, F)$ is represented by the de Rham complex of $F$:
\[\dR(F):=0\ra \Omega^0_X\otimes F\xrightarrow{d_F}
\Omega^1_X\otimes F\xrightarrow{d_F} \cdots\] 
where the differentials $d_F: \Omega^j_X\otimes F \ra
\Omega^{j+1}_X\otimes F$ are given by 
\[d_F(\omega\otimes f)=d(\omega)\otimes f +(-1)^j \omega \wedge
\nabla(f).\] 
In particular, if $F=\cO_X$, the de Rham complex $\dR(\cO_X)$ is the (algebraic) de Rham complex of $X$.

\paragraph
\label{center_D}
The following general properties of the ring of crystalline
differential operators are well known, see~\cite{BezMirRum}.  {\em A
  priori} $D_X$ is constructed as a sheaf of (left) $\cO_X$-modules on
$X$.  However, the center of $D_X$ can be canonically identified with
the structure sheaf of the cotangent bundle $T^*X'$ of $X'$, the
Frobenius twist of $X$.  We review this identification.

Since the map $\pi: X'\ra X$ is an isomorphism, it
induces an isomorphism of sheaves $\pi_*:T_{X'}\ra T_X$, which in local
coordinates can be described as 
\[ \pi_*\left(\frac{\del}{\del (x^p)}\right)=\frac{\del}{\del x}. \]
(Obviously, the map is not $k$-linear.)
For any vector field $\del\in T_X$ its $p$-th power $\del^p\in D_X$
acts again as a derivation on the functions on $X$.  Therefore it can
be regarded as another vector field, denoted by $\del^{[p]}$.  In Lie-theoretic terms, 
vector fields on $X$ form a $p$-restricted Lie algebra.

Consider the map $\iota:T_{X'} \ra D_X$ obtained as the composite
\[ T_{X'} \xrightarrow{\pi_*} T_X \xrightarrow{\del \mapsto
  \del^p-\del^{[p]}} D_X. \]
It is well known~\cite[1.3.1]{BezMirRum} that $\iota$ is
$\cO_{X'}$-linear, and the image is in the center of $D_X$. Since
$f^p$ lies in the center of $D_X$ for $f\in \OO_X$, $\iota$ extends to
a homomorphism $\bar{\iota}:\cO_{T^*X'}\ra Z(D_X)$. One can
check that this map is an isomorphism~\cite[1.3.2]{BezMirRum}.

\paragraph\label{D and DX}
As a consequence, we can regard $D_X$ as a sheaf of algebras on
$T^*X'$ instead of on $X$. To avoid confusion we will denote
$D_X$, when regarded as a sheaf on $T^*X'$, by $D$.  The
projection to the zero section $\phi:T^*X'\ra X'$ identifies $D$
and $D_X$; we have $\phi_*D=F_* D_X$.

\paragraph
\label{Morita}
The algebra $D_X$ is an Azumaya algebra over its center, and so $D$ is a sheaf of 
Azumaya algebras over $T^*X'$.  The restriction
of $D$ to the zero section, $D|_{X'}=i^*D$, is a split Azumaya
algebra, with splitting module $F_*\OO_X$. In other words
$D|_{X'}=\sEnd_{X'}(F_*\OO_X)$. Therefore, $\OO_{X'}$ and $D|_{X'}$
are Morita equivalent: the functors
\begin{align*}
m_* &: \gMod\mbox{-}\OO_{X'}\rightarrow \gMod\mbox{-}D|_{X'}, \quad\quad M\mapsto M\otimes
F_*\OO_X \\
\intertext{and}
m^*& :\gMod\mbox{-}D|_{X'}\rightarrow \gMod\mbox{-}\OO_{X'}, \quad \quad N\mapsto
\sHom_{D|_{X'}}(F_*\OO_X,N)
\end{align*}
are inverse to one another. These functors give rise to an equivalence
between the corresponding derived categories $\D(X',
D_{X'})$ and $\D(X')$.  We shall denote the induced functors on
derived categories by $m_*$ and $m^*$ as well.

\section{Azumaya schemes}
\label{sec:Azumaya}

Our results are naturally formulated in the framework of schemes
equipped with Azumaya algebras. Such objects are called `twisted
spaces' in \cite{CalThesis}, but we prefer the more descriptive name
\emph{Azumaya schemes}. It is worth noting that one can pass from
Azumaya schemes to the corresponding $\Gm$-gerbes (see Remark~\ref{rem:gerbes}); the
results below can be easily restated in the language of gerbes. One
advantage of Azumaya algebras is that they are somewhat more explicit
than gerbes.  

\paragraph{\bf Definition.} An \emph{Azumaya scheme} is a pair $(S,
\cA)$ where $S$ is a scheme and $\cA$ is a (coherent) sheaf of
Azumaya algebra over $S$. Given a property of schemes that is local in
the \'etale topology we say that $(S,\cA)$ has this property if and
only if $S$ does.

\paragraph 
Suppose $(S,\cA)$ is an Azumaya scheme and $f:S'\to S$ is a morphism
from another scheme $S'$ to $S$. Then the pullback $f^*\cA$ is an
Azumaya algebra over $S'$.  We could consider the category whose objects
are Azumaya schemes $(S,\cA)$, and where morphisms $(S',\cA')\to
(S,\cA)$ are pairs
\[(f:S'\to S,\phi:f^*\cA\simeq\cA').\] 
However, such `strict' category of Azumaya schemes is too restrictive
from the view-point of Azumaya algebras: one should work with Morita
equivalences instead of isomorphisms. This leads to the following
definition.

\paragraph{\bf Definition.} 
A morphism (or, more precisely, a $1$-morphism) 
\[ (S',\cA')\to (S,\cA) \]
of Azumaya schemes is a pair $(f,E)$, where $f:S'\to S$ is a morphism
of schemes and $E$ is an $f^*(\cA)^{op}\otimes\cA'$-module that
provides a Morita equivalence between the Azumaya algebras $\cA'$ and
$f^*\cA$.

Given two $1$-morphisms $(f_1,E_1)$ and $(f_2,E_2)$ between Azumaya
schemes $(S',\cA')$ and $(S,\cA)$, $2$-morphisms between them exist
only if $f_1=f_2$. In the case $f_1=f_2$, such $2$-morphisms are
isomorphisms between $f_1^*(\cA)^{op}\otimes\cA'$-modules $E_1$ and
$E_2$.

\paragraph{\bf Remark.} 
The condition that $E$ provides a Morita equivalence between the
Azumaya algebras $\cA'$ and $f^*\cA$ implies that $E$ is a vector
bundle on $S'$ whose rank at $y'\in S$ is the product of the ranks of
$\cA'$ at $y'$ and $\cA$ at $f(y')$.  (The rank of an Azumaya algebra
$\cA$ is $n$ if its rank as an $\cO_S$-module is $n^2$.)

\paragraph 
It is easy to see that Azumaya schemes form a $2$-category with
respect to the natural composition of $1$-morphisms. The category is
equipped with a forgetful functor to the category of schemes, which
sends the Azumaya scheme $(S,\cA)$ to $S$. We think of the Azumaya
scheme $(S,\cA)$ as a `twisting' of the scheme $S$; accordingly, we
call the forgetful functor the \emph{untwisting}.

Conversely, for any scheme $S$ we have the `untwisted' Azumaya scheme
$(S,\cO_S)$; if no confusion is likely we still denote this Azumaya
scheme by $S$. Note however that the corresponding functor from the
category of schemes to the $2$-category of Azumaya schemes is not
fully faithful (see Example~\ref{ex:twist by a line bundle}).

\paragraph 
From now on we suppose that all our Azumaya schemes are of finite
type over a fixed ground field; this ensures that direct/inverse image
functors are defined and well-behaved.  A \emph{quasicoherent sheaf}
on an Azumaya scheme $(S,\cA)$ is a sheaf of $\cA$-modules that is
quasi-coherent as a $\cO_S$-module. Quasicoherent sheaves on $(S,\cA)$
form an abelian category with enough injectives; denote the
corresponding derived category by $\D(S,\cA)$. Of course, if
$(S,\cO_S)$ is the untwisted Azumaya scheme, the category $\D(S,\cO_S)$
is the usual derived category $\D(S)$.

Let $\phi=(f,E)$ be a morphism of Azumaya schemes $(S',\cA')\to
(S,\cA)$. The direct image functor $f_*:\D(S')\to\D(S)$ and the
inverse image functor $f^*:\D(S)\to\D(S')$ are naturally upgraded to
functors $f_*:\D(S',f^*\cA)\to\D(S,\cA)$ and
$f^*:\D(S,\cA)\to\D(S',\cA')$; composing them with the Morita
equivalence
\[\D(S',\cA')\simeq\D(S',f^*\cA)\]
given by $E$, we obtain the direct image functor
\begin{align*}\phi_*:\D(S',\cA')&\to\D(S,\cA)\\M'
&\mapsto f_*(E^\vee\otimes_{\cA'} M')\\
\intertext{and the inverse image functor}
\phi^*:\D(S,\cA)&\to\D(S',\cA')
\\M&\mapsto E\otimes_{f^*\cA} f^*(M).
\end{align*}
It is easy to see that $\phi^*$ is the left adjoint of $\phi_*$.

\paragraph 
Let us now consider the `extraordinary' inverse image functor. Let
$\phi$ be as above. Generally speaking, the inverse image functor
$f^!$ is defined only on bounded below categories
\[f^!:\D^+(S)\to\D^+(S');\]
accordingly, we obtain the inverse image functor 
\begin{align*}
\phi^!:\D^+(S,\cA)&\to\D^+(S',\cA')
\\M&\mapsto E\otimes_{f^*\cA} f^!(M)
\end{align*}
between the bounded below categories of the Azumaya schemes. If $f$ is
Gorenstein (or, more generally, if the relative dualizing complex of
$S'$ over $S$ is perfect), the functor $f^!$ can be extended to the
unbounded categories, and then $\phi^!$ also extends to the inverse
image functor
\[\phi^!:\D(S,\cA)\to\D(S',\cA').\]
For instance, this is always the case for morphisms between smooth
Azumaya schemes.

If the map $\phi$ is proper (that is, if $f$ is proper), then $\phi^!$ is the right adjoint of $\phi_*$.

\paragraph
\label{ex:twist by a line bundle} 
{\bf Example.} Let $(S,\cA)$ be an Azumaya scheme. The identity
$1$-morph\-ism of $(S,\cA)$ is the pair $(id_S,\cA)$, where we view
$\cA$ as a $\cA^{op}\otimes\cA$-module in the natural way. If $\sL$
is a line bundle on the (untwisted) scheme $S$, the pair
\[\phi_\sL=(id_S,\sL\otimes\cA)\] 
is a $1$-automorphism of $(S,\cA)$ as well; in this way, we obtain all
$1$-automor\-phisms of $(S,\cA)$ that act trivially on its `untwisting'
$S$. In particular, even if $\cA=\cO_S$, and $(S,\cA)$ is the usual
`untwisted' scheme, it acquires some `twisted' automorphisms
($\phi_\sL$) when considered as an Azumaya scheme.

The pull-back and push-forward functors for the morphism $\phi_\sL$
are then given by
\begin{align*}
\phi^*=\phi^!:\D(S,\cA)&\to\D(S,\cA)&\phi_*:\D(S,\cA)&\to\D(S,\cA)\\
M&\mapsto M\otimes\sL&M&\mapsto M\otimes\sL^{-1}.
\end{align*}

 \paragraph 
\label{sec:fiber product}
Suppose we are given three Azumaya schemes $\bX = (X,\cA_X)$,
$\bY = (Y,\cA_Y)$, and $\bS = (S,\cB)$, and morphisms 
\[ \phi_X = (f_X, E_X): \bX \ra \bS, \quad \phi_Y = (f_Y, E_Y):\bY\ra
\bS. \] 
The fibered product $W:=X\times_S Y$ is equipped with a natural
Azumaya algebra $\cA$, making it into an Azumaya space $\bW = (W,
\cA)$. Explicitly,
\[p_1^*(E_X)\otimes p_2^*(E_Y)\]
is a module over 
\[p_1^*(\cA_X)\otimes p_2^*(\cA_Y)\otimes(p^*(\cB)^{\otimes2})^{op},\]
and we let $\cA$ be the algebra of its endomorphisms as a
$(p^*(\cB)^{\otimes2})^{op}$-module. Here $p_1$, $p_2$, and $p$ are
the projections of $W$ onto $X$, $Y$, and $S$, respectively.

The four Azumaya algebras $\cA$, $p_1^*(\cA_X)$, $p_2^*(\cA_Y)$, and
$p^*(\cB)$ are naturally Morita equivalent; this gives a commutative
square of Azumaya schemes
\[\xymatrix{\bW \ar[d] \ar[r] & \bX \ar[d]\\
\bY\ar[r]& \bS.}\]
It is not hard to see that this square is $1$-Cartesian, so that
$\bW$ is the fiber product of $\bX$ and $\bY$ over
$\bS$.

A particular case of this construction is when $\cA_X =  f_X^* \cB$,
$\cA_Y = f_Y^*\cB$.  In this case the algebra $\cA$ can naturally be
identified with $p^*\cB$.

\paragraph{\bf Remark.} 
Note that the $2$-category of Azumaya schemes does not have a final
object. For this reason, while fiber products exist, Cartesian
products do not.

\paragraph 
\label{sec:intersection}
Suppose now in addition that the two Azumaya schemes $\bX$, $\bY$ are
untwisted, so that $\cA_X=\cO_X$, $\cA_Y = \cO_Y$.  In this case, the
Azumaya scheme $(W,\cA)$ is also split. However, it carries two
natural splittings: one of them compatible with the splitting of the
Azumaya scheme $\bX$ and the other compatible with the
splitting of $\bY$. The difference between the two splittings
is a line bundle $\sL$ on $S$. Explicitly,
\[\sL=\sHom_{p^*(\cB)^{op}}(p_1^*(E_1),p_2^*(E_2))\]
measures the discrepancy between the splitting of the pullbacks
$f_X^*(\cB)$ and $f_Y^*(\cB)$.

In this way, we obtain a commutative diagram
\[\xymatrix{
&(W,\cO_W)\ar[d]^{\sim}\ar[rd]&\\
(W,\cO_W)\ar[ru]_{\sim}^{\phi_\sL}\ar[rd]\ar_{\sim}[r]&(W,\cA)\ar[r]\ar[d]&(X,\cO_X)\ar[d]\\
&(Y,\cO_Y)\ar[r]&(S,\cB)} \] 
of Azumaya schemes. Here the morphism $\phi_\sL$ is constructed in
Example~\ref{ex:twist by a line bundle}.

\paragraph 
All of the above constructions admit differential-graded analogues. In
particular in the situation of Section~\ref{sec:fiber product} we
may consider the derived fiber product of the Azumaya
schemes. Specifically, the derived fiber product of $X$ and $Y$ over
$S$ is naturally equipped with a sheaf of Azumaya algebras, and
we thus obtain an \emph{Azumaya dg-scheme} $W^R$.

This is the only way in which Azumaya dg-schemes appear in this
paper. In fact, for our purposes the most important situation is as
in~(\ref{sec:intersection}): the Azumaya schemes $\bX$ and $\bY$
are non-twisted, and thus their derived fiber product $W^R$ is split (that
is, $1$-isomorphic to an untwisted dg-scheme); however, there are two
splittings that differ by a line bundle on the dg-scheme $W^R$. For
details, see~(\ref{subsec:splitting})-(\ref{subsec:lineb}).

\paragraph
\label{rem:gerbes}
{\bf Remark.} It is well known that an Azumaya algebra on a scheme $S$
yields a class in the Brauer group $H^2(S,\cO_S^\times)$, which can be
interpreted geometrically as a $\Gm$-gerbe on $S$. (Here $\Gm$ is the
multiplicative group.)

Explicitly, to an Azumaya scheme $(S,\cA)$, we assign the stack
$\Split(\cA)$ such that morphisms $T\to\Split(\cA)$ from a test scheme
$T$ are pairs $(f,E)$, where $f:T\to S$, and $E$ is the splitting of
the Azumaya algebra $f^*\cA$. Equivalently, $(f,E^\vee)$ is a morphism
of Azumaya schemes $(T,\cO_T)\to(S,\cA)$.

It is clear that $\Split(\cA)$ is equipped with a natural morphism to
$S$, and that this morphism turns $\Split(\cA)$ into a $\Gm$-gerbe
over $S$. In this way, we obtain a fully faithful embedding of the
$2$-category of Azumaya schemes into the $2$-category of $\Gm$-gerbes.

Since $\Split(\cA)$ is a $\Gm$-gerbe, the category $\D(\Split(\cA))$
decomposes into direct summands indexed by the characters of $\Gm$.
The summand corresponding to the tautological character is identified
with $\D(S,\cA)$. This identification agrees with the pullback and
pushforward functors for a morphism between Azumaya schemes.

\section{Twisted derived intersections}
\label{sec:twiderint}

In this section we discuss derived intersections in twisted spaces,
generalizing the results of~\cite{AriCal} and~\cite{AriCalHab}.  We
argue that there is a natural twisted HKR class with similar
properties to those of the usual one.  The most important observation
is that, with certain assumptions, the twisted derived intersection
differs only by a line bundle from the untwisted one.   
\medskip

\noindent
We begin our discussion with the general framework for the derived
intersection problems we study.

\paragraph
\label{subsec:context}
Let $S$ be a smooth variety, and let $X$ and $Y$ be smooth subvarieties
of $S$, with corresponding closed embeddings $i$ and $j$.  We shall
denote by $W^0$ and $W$ their underived and derived intersections,
\[ W^0 = X\times_S Y, \quad\quad W = X\times_S^R Y. \] 
We assume given an Azumaya algebra $\cA$ on $S$ which allows us to
construct the Azumaya space $\bS = (S, \cA)$.  Let $\bX$, $\bY$ be the
Azumaya schemes $(X, \cA|_X)$, $(Y, \cA|_Y)$.  They are naturally
smooth subvarieties of $\bS$ with corresponding closed embeddings by
$\bi$, $\bj$.  The derived intersection $\bW$ of $\bX$ and $\bY$
inside $\bS$ is given by 
\[ \bW = \bX \times_\bS^R \bY = (W, \cA|_W), \] 
see~(\ref{sec:fiber product}).  

\paragraph
\label{subsec:splitting}
Now assume that the restrictions of $\cA$ to $X$ and to $Y$ are split
by vector bundles $E_X$ and $E_Y$ on $X$ and $Y$, respectively:
\[ \cA|_X \iso \End_X(E_X), \quad\quad \cA|_Y \iso \End_Y(E_Y).\] 
This data gives rise to isomorphisms of Azumaya schemes $m:X
\xrightarrow{\sim} \bX$, $n:Y\xrightarrow{\sim} \bY$ which allow us to
regard $X$ and $Y$ as smooth subvarieties of either the untwisted
space $S = (S, \cO_S)$ or of the twisted space $\bS = (S, \cA)$.  In
the latter case the closed embedding $i':X \hookrightarrow \bS$ is
given by the composite $X \xrightarrow{m} \bX \xrightarrow{\bi} \bS$.
A similar composition gives the embedding $j'$ of $Y$ into $\bS$.

\paragraph
The given splittings of $\cA|_X$ and $\cA|_Y$ give rise to splittings
of $\cA|_W$.  These can be regarded as isomorphisms $f$ and $g$
between $W$ and $\bW$.  The spaces and morphisms we have
discussed so far fit in the following diagram:
\[\xymatrix{W \ar@/_/[dr]_{g} \ar@/^/[dr]^{f} \ar[dd]_{q}
\ar[rr]^p& & X\ar[d]^m \ar@/^2pc/[dd]^{i'}\\ & \bW
\ar[r]^{\bar{p}}\ar[d]^{\bar{q}} & \bX\ar[d]^{\bar{i}}\\
Y\ar@/_2pc/[rr]^{j'}\ar[r]^n& \bY\ar[r]^{\bar{j}} & \bS.\\} \]
Note that we are not assuming that $E_X|_W \iso E_Y|_W$, hence the
morphisms $f$ and $g$ are not necessarily isomorphic.  In particular,
the entire diagram is {\em not} commutative, even though the
right-upper and left-lower trapezoids are.

\paragraph
We are interested in comparing the derived intersection $W$ of $X$ and
$Y$ inside $S$ with its twisted counterpart, when we consider the
intersection as taking place inside $\bS$.  Denote by $W'$ this latter 
dg scheme,
\[ W' = X \times_\bS^R Y. \] 
The structure sheaves of $W$ and $W'$ (regarded as dg schemes with
morphisms to $X$ and to $Y$) are the kernels of $j^*i_*$ and
$j'{}^*i'_*$, respectively, both regarded as functors $\D(X) \ra
\D(Y)$.  (For a general discussion of the point of view that structure
complexes of derived intersections should be understood as kernels of
functors see~\cite[Section 2]{AriCalHab}.)

Note that we would normally think of $W$, $W'$ as dg schemes over
$X\times Y$ because they are both endowed with morphisms to $X$ and
$Y$.  However, since the category of Azumaya spaces does not have
absolute products (as opposed to fibered products), it is more natural
to think of them as dg schemes over the pair $(X, Y)$: such a spaces
is endowed with two morphisms, one to $X$ and one to $Y$, and
morphisms in this category must make the obvious diagrams commutative.

\paragraph
\label{subsec:lineb}
Regarding derived intersections as kernels makes it clear that the dg
scheme $W'$ can be taken to be $W$ as an abstract dg scheme, however
one of the two morphisms $p$ or $q$ has to be modified in order to fix
the commutativity of the resulting diagram.  Let $\tau$ be the
automorphism (in the category of Azumaya dg schemes) of $W$ given by
$\tau = g^{-1} \circ f$.  As an abstract Azumaya dg scheme, $W'$ is
isomorphic to $W$.  However as a space over $(X, Y)$ it is different
from $W$ -- its projection maps to $X$ and to $Y$ are given by $p$ and
$q\circ\tau$, respectively, instead of $p$ and $q$.  Note that the
resulting diagram commutes:
\[  \xymatrix{ W=W' \ar[dr]^{f} \ar[dd]_{q\circ \tau}
\ar[rr]^p& & X\ar[d]^m \ar@/^2pc/[dd]^{i'}\\ & \bW
\ar[r]^{\bar{p}}\ar[d]^{\bar{q}} & \bX\ar[d]^{\bar{i}}\\
Y\ar@/_2pc/[rr]^{j'}\ar[r]^n& \bY\ar[r]^{\bar{j}} & \bS.\\} \]
It is easy to understand the automorphism $\tau$ of $W$.  At the
level of underlying dg schemes $\tau$ is the identity, hence $\tau$ must be
given by a line bundle $\sL$ on $W$, see~(\ref{ex:twist by a line bundle}).
Rephrasing the commutativity of the above diagram we get
\[ j'{}^*i'_* (\,-\,) \iso q_* (p^*(\,-\,) \otimes \sL). \]
We will call $\sL$ the {\em associated line bundle} of the twisted
derived intersection problem.

\paragraph {\bf Remark.}
Line bundles on a dg scheme $W$ are classified, as in the
classical case, by 
\[ Pic(W) = H^1(W, \cO_W^\times) = H^1(W^0, \cO_W^\times). \]
If $\cO_W$ is formal and concentrated in non-positive degrees then the
complex $\cO_W^\times$ of invertible elements in $\cO_W$ consists of
those whose degree zero part is invertible.  (Indeed, $\cO_W$ is a
nilpotent algebra over $\cO_{W^0}$.)  We then
have
\[ \Pic(W) = \Pic(W^0) \oplus \bigoplus_{i\geq 1} H^{i+1}(W^0,
\coH^{-i}(\cO_W) ), \]
where $\coH^{-i}(\cO_W)$ is the $-i$-th cohomology sheaf of $\cO_W$.

In particular $\Pic(W^0)$ and $H^2(W^0, \coH^{-1}(\cO_W))$ are direct
summands in $\Pic(W)$.  It is easy to understand these groups
geometrically. If $\sL$ is a line bundle on $W$, its zeroth cohomology
sheaf $\coH^0(\sL)$ is a line bundle on $W^0$, and this gives the map
$\Pic(W) \ra \Pic(W^0)$.  (One can also understand this map as the
pull-back via the embedding $W^0 \hookrightarrow W$.)  

Since $\cO_W$ is formal there is a projection map $\pi:W\ra W^0$.
Even though $\pi_*\cO_W$ is a formal complex on $W_0$, for a line
bundle $\sL$ on $W$ the complex $\pi_*\sL$ need not be formal.  The
first obstruction to the formality of $\pi_*\sL$ is its HKR class
(see~(\ref{hkr:class}) below), the obstruction to the splitting of
$\tau^{\geq -1}(\pi_*\sL)$, which is a class in
\[ Ext^2_{W^0}(\coH^0(\pi_*\sL), \coH^{-1}(\pi_*\sL)) = H^2(W^0,
\coH^{-1}(\cO_W)). \] 
In this way we realize the HKR class for a line bundle $\sL$ as a
component of the class of $\sL$ in $\Pic(W)$.  \medskip

\noindent
The main result of this section is the following.

\begin{Theorem}
\label{thm:twiderint}
The associated line bundle $\sL$ is trivial (or, equivalently, $W\iso
W'$ as Azumaya dg schemes over $X\times Y$) if conditions (a)--(c)
below are satisfied
\begin{itemize}
\item[(a)] $E_X|_{W^0} \iso E_Y|_{W^0}$;
\item[(b)] the map $i'$ is split to first order; and
\item[(c)] the natural map $N_{W^0/Y} \ra N_{X/S}|_{W^0}$ of vector bundles
  on $W^0$ is split.
\end{itemize}
If this is the case then both $W$ and $W'$ are formal dg schemes over
$(X, Y)$, $W \iso W' \iso \Tot_{W^0}(\bbE[-1])$, where 
\[ \bbE = \frac{T_S}{T_X + T_Y} \]
is the excess intersection bundle for the original intersection
problem $X\cap Y$ inside $S$.

Moreover, if $X = Y$, $i'=j'$ (i.e., $i = j$ and $E_X = E_Y$) and we
assume that $i$ is split to first order, then the converse also holds: $W
\iso W'$ if and only if $i'$ is split to first order.
\end{Theorem}

\begin{Proof}
The proof is essentially the same as the one of the corresponding
result in~\cite{AriCal},~\cite{AriCalHab}.  Therefore we only review
the main points of the proof highlighting the aspects that are special
to the twisted case.

Let us first review the proof in the untwisted case.  There the proof
shows that $W \iso \Tot_{W^0}(\bbE[-1])$ if $i$ is split to first
order and condition (c) is satisfied. We shall describe the proof in
geometric language.  We consider the sequence of maps of dg spaces over
$(X, Y)$
\begin{align*}
\Tot_{W^0}(\bbE[-1]) & \ra \Tot_{W^0}(N_{X/S}|_{W^0}[-1]) \iso
\Tot_X(N_{X/S}[-1]) \times_X W^0 \iso \\
&\iso (X\times_S X) \times_X W^0 \iso (X\times_S Y)\times_Y W^0 \ra
X\times_S Y. 
\end{align*}
Here the first map comes from the splitting of the sequence on $W^0$
\[ 0 \ra N_{W^0/Y} \ra N_{X/S}|_{W^0} \ra \bbE \ra 0, \]
while the second isomorphism is the main result of~\cite{AriCal}, and
its existence follows from the fact that $i$ is split to first order.
(Throughout this proof we will implicitly assume that all fiber
products are derived, dropping the superscript ``$R$'' for simplicity.)

Having constructed a map $\Tot_{W^0}(\bbE[-1]) \ra X\times_S Y$ of dg spaces over
$(X,Y)$, checking that it is a quasi-isomorphism is a local question
which can be settled using Koszul resolutions to complete the proof.

Moving over to the twisted case, let $\cA$ be an Azumaya algebra over
$S$.  The same exact proof shows that conditions (b) and (c) imply
that $\bW \iso \Tot_{\bW^0}(\bbE[-1])$ as $(\bX, \bY)$ spaces, where
placing a bar over an $S$-space $X$ means endowing it with the Azumaya
algebra $\cA|_X$.  Note that the condition that $\bi$ splits to first
order means that the splitting $\pi: X^{(1)} \ra X$ must satisfy the
condition that $\pi^*(\cA|_X) \iso \cA|_{X^{(1)}}$.

The entire discussion above makes no assumptions on the splitting of
$\cA$ on $X$ or $Y$.  Now assume furthermore that we have fixed
splittings $E_X$, $E_Y$ of $\cA|_X$, $\cA|_Y$ as
in~(\ref{subsec:splitting}).   Condition (a) of the theorem means
that there is an isomorphism $W^0 \iso \bW^0$ making a commutative
diagram of Azumaya spaces 
\[ \xymatrix{ W^0 \ar[dr] \ar[dd] \ar[rr]& & X\ar[d]\\ & \bW^0
  \ar[r]\ar[d] & \bX\ar[d]\\
  Y\ar[r]& \bY\ar[r] & \bS.\\} \] 
(All the spaces in this diagram are non-derived.)  Thus $W^0$, as an
$(X, Y)$-space, is isomorphic to $\bW^0$ as an $(\bX, \bY)$-space.
The previous discussion shows that conditions (b) and (c) imply that
$\bW$ is the total space of $\bbE[-1]$ over $\bW^0$, as a space over
$(\bX, \bY)$.  Thus $W'$, as the pull-back of $\bW$ over the map of
pairs of spaces $(X, Y) \ra (\bX, \bY)$, is isomorphic to
$\Tot_{W^0}(\bbE[-1])$ as an $(X, Y)$-space.  Since the untwisted map
$i$ is split and condition (c) holds we know that the untwisted
intersection $W$ is also isomorphic to $\Tot_{W^0}(\bbE[-1])$.  This
completes the proof of the fact that 
\[ W \iso W' \iso \Tot_{W^0}(\bbE[-1]). \]

The reverse implication (the second part of the theorem) follows from
the discussion below of twisted HKR classes, see Theorem~\ref{twisted}. \qed
\end{Proof}

\paragraph
\label{hkr:class}
In the course of the above proof we used a generalization of
Theorem~\ref{AriCal} which applies for Azumaya spaces.  We discuss
this generalization below.

Consider the following set up.  Let $i:X\inj S$ be a closed
embedding of smooth schemes over a field $\bbk$ of characteristic either
zero or greater than the codimension of $X$ in $S$. Denote
the first infinitesimal neighborhood of $X$ in $S$ by $X^{(1)}$.  

Let $\sA$ be an Azumaya algebra over $S$, providing a map
$\bi:(X,\sA|_X)\inj (S,\sA)$ of Azumaya schemes. Let $E$ be a
coherent sheaf on $(X, \sA|_X)$ that is locally free as an $\cO_X$-module (a `twisted vector bundle' for the
twisting $\sA|_X$).  We are interested in the object
$\bi^*\bi_*E$ of $\D(X,\sA|_X)$. A local computation with the Koszul complex
(similar to the one in \cite{CalKatSha}) shows that its cohomology sheaves
are given by
\[ \coH^{-j}(\bi^*\bi_*E) = E \otimes \wedge^j N^\chk, \]
where $N^{\chk}=N^\chk_{X/S}$ denotes the conormal bundle of $X$
in $S$. Therefore we have a triangle in $\D(X,\sA|_X)$
$$E\otimes N^{\chk}[1]\ra \tau^{\ge -1}\bi^*\bi_*E\ra E\ra E\otimes N^{\chk}[2].$$
We denote the rightmost map in this triangle by
\[\alpha_E\in\Ext^2_{\sA|_X}(E,E\otimes N^\chk)\]
and call it the \emph{(twisted) HKR class} of $E$, in keeping with the notation
of~\cite{AriCal}.
\smallskip

\noindent
The HKR class satisfies a number of simple properties.

\begin{Proposition}\label{propHKR}
\begin{itemize}
\item[(1)] Suppose $S_2$ is another smooth scheme, with a smooth
  morphism $f:S_2\to S$.  Let $X_2 = f^{-1} X$.  Consider the Azumaya
  algebra $\cA_2:=f^*\cA$ on $S_2$. Set $E_2:=f^*(E)$; it is a
  coherent sheaf on $(X_2,\cA_2|_{X_2})$ that is locally free over
  $\cO_{X_2}$. Then
\[\alpha_{E_2}=f^*(\alpha_E)\in\Ext^2_{{\sA_2}|_{X_2}}(E_2,E_2\otimes N_2^\chk),\]
where $N_2=f^*(N)$ is the normal bundle to $X_2\subset S_2$.
\item[(2)] Suppose that $E_1,E_2$ are two coherent sheaves on
  $(X,\cA|_X)$ that are locally free as $\cO_X$-modules. Then the HKR
  class of $E_1\oplus E_2$ is the image of
\[(\alpha_{E_1},\alpha_{E_2})\in\Ext^2_{\sA|_X}(E_1,E_1\otimes N^\chk)\oplus \Ext^2_{\sA|_X}(E_2,E_2\otimes N^\chk)\]
under the natural map to $\Ext^2(E_1\oplus E_2,(E_1\oplus E_2)\otimes N^\chk)$.
\item[(3)] Suppose $\cA_1$ and $\cA_2$ are two Azumaya algebras on
  $S$, and let $E_i$ be a coherent sheaf on $(X,\cA_i|_X)$ that is
  locally free as an $\cO_X$ module ($i=1,2$). The HKR class of the
  coherent sheaf $E_1\otimes E_2$ on $(X,(\cA_1\otimes\cA_2)|_X)$ is
  given by
\[\alpha_{E_1\otimes E_2}=\alpha_{E_1}\otimes id_{E_2}+ id_{E_2}\otimes\alpha_{E_1}.\]

\item[(4)] In the situation of the previous statement, the HKR class
  of the coherent sheaf $\sHom(E_1,E_2)$ on
  $(X,(\cA_1^{op}\otimes\cA_2)|_X)$ is given by
\[\alpha_{\sHom(E_1,E_2)}=-\alpha_{E_1}\otimes id_{E_2}+ id_{E_2}\otimes\alpha_{E_1}.\]
Here $\cA_1^{op}$ is the opposite of the Azumaya algebra $\cA_1$. \qed
\end{itemize}
\end{Proposition}

\paragraph While we defined the HKR class using the object $\bi^*\bi_*E$, we can use 
$\bi^!\bi_*E$ instead. Indeed, local computation with the Koszul complex
shows that 
\[ \coH^{j}(\bi^*\bi_*E) = E \otimes \wedge^j N, \]
Therefore we have a triangle in $\D(X,\sA|_X)$
$$E\ra \tau^{\le 1}\bi^*\bi_*E\ra E\otimes N[-1]\ra E[1].$$
\begin{Lemma}
The rightmost map in this triangle is equal to $\alpha_E$.
\end{Lemma}

\begin{Proof} Consider $E^\vee=\sHom(E,\cO_X)$ as a module over the
  Azumaya algebra $\sA^{op}|_X$.  If we let $\bi^{op}$ be the map of
  Azumaya spaces $\bi^{op}:(X, \sA^{op}|_X) \ra (S, \sA^{op})$ we have
\[\bi^!\bi_*E=\RsHom(\bi^{op,*}i^{op}_*(E^\vee), \sA^{op});\]
now the claim follows from Proposition~\ref{propHKR}(4).
\end{Proof}

\paragraph{} Finally, we can view the HKR class as the obstruction to extending $E$ to the first infinitesimal neighborhood of the
Azumaya scheme $(X,\cA|_X)$. More precisely, local extensions of $E$ form a gerbe of extensions, which we denote by $\EG$; and
the class of this gerbe $[\EG]$ is equal to $\alpha_E$.
In the `untwisted' case of ordinary schemes, this gerbe was considered in \cite{AriCal}; the case of Azumaya schemes is completely 
analogous. 

Explicitly, we construct the gerbe $\EG$ on $X$ as follows. Given an open subset $U\subset X$, we denote by $U^{(1)}\subset X^{(1)}$ its first infinitesimal neighborhood.
By definition, $\EG(U)$ is the groupoid whose objects are pairs consisting of a $\cA|_{U^{(1)}}$-module $E^{(1)}$ that is
locally free over $\cO_{U^{(1)}}$ and an isomorphism between the restriction $E^{(1)}|_U$ and the $\cA|_U$-module $E|_U$. 
(To simplify notation, we usually omit this isomorphism.) As $U\subset X$ varies, the groupoids $\EG(U)$ form a sheaf of groupoids
$\EG$ on $X$. It is locally non-empty (locally, sections exist) and locally connected (any two sections are locally isomorphic).
Moreover, the sheaf of automorphisms of any section is naturally identified with
\[\sHom_{\sA|_X}(E,E\otimes N^\vee)=N^\vee.\]
Thus, $\EG$ is a gerbe over the sheaf $N^\vee$. The class of this gerbe
\[[\EG]\in H^2(X, N^\vee)\] 
is equal to $\alpha_E$ by the argument that is completely parallel to the proof of \cite[Proposition~2.6]{AriCal}.

\medskip

\begin{Theorem}
\label{twisted}
In the above setting, assume that the normal bundle $N=N_{X/S}$ extends
to the first infinitesimal neighborhood $X^{(1)}$. Then the following statements are equivalent.
\begin{itemize}            
\item[(1)] The twisted HKR class $\alpha_E$ of $E$ vanishes.
\item[(2)] There exists a formality isomorphism
  \[\bi^*\bi_*E\iso E\otimes \S(N^{\chk}[1])\] in
  $\D(X,\sA|_X)$ (where the action of $\sA|_X$ on the right-hand side comes from its action on $E$).
\item[(3)] There exists a formality isomorphism
\[\bi^!\bi_*E\iso E\otimes \S(N[-1])\] in
  $\D(X,\sA|_X)$ (where the action of $\sA|_X$ on the right-hand side comes from its action on $E$).
\item[(4)] The sheaf $E$ extends to the first infinitesimal neighborhood $X^{(1)}$. Recall that by this we mean
that there exists a coherent sheaf $E^{(1)}$ on the Azumaya scheme $(X^{(1)},\sA|_{X^{(1)}})$ that is locally
free over $\cO_{X^{(1)}}$ and whose restriction
to $(X,\sA|_X)$ is isomorphic to $E$. 
\end{itemize}
\end{Theorem}

\begin{Proof} 
This theorem is entirely similar to the original result~\cite[Theorem
0.7]{AriCal}, and we shall leave to the reader the task of filling in
the details of the proof. \qed 
\end{Proof}

\paragraph
In the situation of Theorem~\ref{twisted} suppose that $E$ splits the
Azumaya algebra $\sA|_X$. Recall that this means that $\rk(E)$ is
equal to the rank of the Azumaya algebra $\sA|_X$. Note that in this
case,
\[\alpha_E\in\Ext^2_{\sA|_X}(E,E\otimes N^\chk)=H^2(X,N^\chk).\] 
An extension $E^{(1)}$, if it exists, splits the Azumaya algebra
$\sA|_{X^{(1)}}$. Geometrically, $E$ gives a morphism of Azumaya schemes 
\[(X,\cO_X)\to(S,\sA),\] while $E^{(1)}$ is an extension of this
morphism to $(X^{(1)},\cO_{X^{(1)}})$. We arrive at the following
conclusion.

\begin{Corollary}\label{twisted splitting}
Under the hypotheses of Theorem~\ref{twisted}, suppose $E$ splits the
Azumaya algebra $\sA|_X$. Then the equivalent conditions of
Theorem~\ref{twisted} are also equivalent to the following condition:
\begin{itemize}
\item[(5)] The morphism \[(X,\cO_X)\to(S,\sA)\]
given by $E$ fits into the diagram
\[(X,\cO_X)\to (X^{(1)},\cO_{X^{(1)}})\to(S,\sA)\]
(where the left arrow is the closed embedding $X\hookrightarrow
X^{(1)}$) such that the underlying morphisms of schemes are the closed
embeddings 
\[X\hookrightarrow X^{(1)}\hookrightarrow S. \]
\end{itemize}
\end{Corollary}

\paragraph
{\bf Remark.}
\label{rem:anysplit}
In particular condition $(5)$ implies that the algebra
$\sA|_{X^{(1)}}$ splits.  However, the converse is
generally not true: it is possible that $\sA|_{X^{(1)}}$ splits, but
there is no splitting that restricts to $E$ on $X$.

On the other hand, if every line bundle on $X$ extends to a line
bundle on $X^{(1)}$ (for example, if the map $X\ra X^{(1)}$ is split),
then condition $(5)$ is equivalent to the weaker condition that
$\sA|_{X^{(1)}}$ splits.  Indeed, let $E^{(1)}$ be any splitting
module for $\sA|_{X^{(1)}}$.  Its restriction $E^{(1)}|_X$ is a
splitting module for $\sA|_X$, thus there exists a line bundle $L$ on
$X$ such that
\[ E \iso E^{(1)}|_X \otimes L. \]
If $L^{(1)}$ is a line bundle on $X^{(1)}$ that extends $L$, then 
$E^{(1)}\otimes L^{(1)}$ is a splitting module on $X^{(1)}$ extending $E$.

\bigskip

\noindent
{\bf \large HKR class in the presence of a section}
\bigskip

\noindent
Let us return to the situation of Corollary~\ref{twisted
  splitting}. Thus $i:X\inj S$ is a closed embedding of smooth schemes
over a field $\bbk$ of characteristic either zero or greater than the
codimension of $X$ in $S$, $\sA$ is an Azumaya algebra on $S$, and $E$
is a splitting module for $\sA|_X$. Above, we have constructed the
twisted HKR class of $E$
\[\alpha_E\in\Ext^2_{\sA|_X}(E,E\otimes N^\chk)=H^2(X,N^\chk).\]
Suppose now that $E$ admits a nowhere vanishing section $s\in
H^0(X,E)$. In this case we can write $\alpha_E$ as a product of two
classes. Let us sketch this construction.

\paragraph
A section $s$ gives an embedding of vector bundles on $X$
\[\cO_X\to E:f\mapsto fs.\]
Denote its cokernel by $E'=E/(\cO_Xs)$. We thus have a short exact sequence
\[0\to\cO_X\to E\to E'\to 0.\]
Let $\beta(E,s)\in\Ext^1(E',\cO_X)$ be the corresponding cohomology
class. Note that $\beta(E,s)$ does not depend on the Azumaya algebra
$\cA$ or its action on $E$.  This is the first of the two classes
which decomposes $\alpha_E$.

\paragraph 
The second cohomology class lies in $\Ext^1((E')^\chk,N^\chk)$. It is
the obstruction to the extension of the pair $(E,s)$ to the first
infinitesimal neighborhood $X^{(1)}$ of $X$. We can describe it as
follows. Since the section $s$ has no zeros, it generates $E$ as an
$\cA|_X$-module. The action of $\cA|_X$ on $s$ gives a homomorphism of
$\cA|_X$-modules $\cA|_X\to E.$ The kernel of this homomorphism is a
left ideal in $\cA|_X$ (the annihilator of $s$), which we denote by
$\cI\subset\cA|_X$.

Consider the problem of deforming $\cI$ to an ideal
$\cI^{(1)}\subset\cA|_{X^{(1)}}$. Locally in the \'etale topology,
$\cA$ splits and such a deformation always exists. Moreover, any two
such deformations differ by a section of
\[\sHom_{\cA|_X}(\cI,(\cA|_X/\cI)\otimes N^\chk).\]
Thus these deformations form a torsor over the sheaf
\[\sHom_{\cA|_X}(\cI,(\cA|_X/\cI)\otimes N^\chk).\]
Note that the Morita equivalence corresponding to the splitting $E$
identifies this sheaf with the sheaf
\[\sHom_{\cO_X}((E')^\chk,N^\chk),\]
so the class of the torsor is an element
\[\gamma(E,s)\in\Ext^1((E')^\chk,N^\chk),\]
as claimed.

\begin{Proposition} The class $\alpha_E\in
  H^2(X,N^\chk)=\Ext^2(N,\cO_X)$ decomposes as the product of
  $\beta(E,s)\in\Ext^1(E',\cO_X)$ and
  $\gamma(E,s)\in\Ext^1((E')^\chk,N^\chk)=\Ext^1(N,E')$.
\end{Proposition}

\begin{Proof} The proof is naturally formulated in the language of gerbes. Let us sketch the argument.

Recall that the HKR class $\alpha_E$ is equal to the class of the gerbe $\EG$ of extensions 
of $E$ to the first infinitesimal neighborhood of the Azumaya scheme $(X,\cA|_X)$.
Similarly, we can represent $\gamma(E,s)$ as the class of a gerbe over the complex
\[C^\bullet=(C^0\to C^1):=(N^\vee\to E\otimes N^\vee).\]
Recall that a gerbe over a complex of sheaves of abelian groups
\[C^\bullet=(C^0\to C^1)\]
is a $C^0$-gerbe $\cG$ and a trivialization of the induced $C^1$-gerbe. Such a trivialization can be given by 
a morphism of sheaves of groupoids $\cG\to\Tors(C^1)$ that is compatible with the $2$-action of $C^0$. Here
$\Tors(C^1)$ is the trivial $C^1$-gerbe: its sections are $C^1$-torsors.

In our situation, consider the $N^\vee$-gerbe $\EG$. We construct a morphism $\EG\to\Tors(E\otimes N^\vee)$ by
sending a local section $E^{(1)}\in\EG(U)$ to the torsor of liftings of $s\in H^0(X,E|_U)$ to a section $s^{(1)}$
of $E^{(1)}$. In this way, we obtain a gerbe $\EG_s$ over the complex $C^\bullet$. The corresponding class
$[\EG_s]$ belongs to the hypercohomology group $H^2(X,C^\bullet)$.
Using the quasi-isomorphism $C^\bullet\to E'\otimes N^\vee[-1]$, we obtain an identification
\[H^2(X,C^\bullet)=H^1(X,E'\otimes N^\vee),\]
and it is easy to see that $[\EG_s]$ corresponds to $\gamma(E,s)$ under this isomorphism.

On the other hand, $\EG$ can be viewed as the gerbe induced by $\EG_s$ under the morphism
\[C^\bullet\to N^\vee.\]
Accordingly, $[\EG]$ is the image of $[\EG_s]$ under this morphism. Now it remains to notice that the composition
\[E'\otimes N^\vee[-1]\simeq C^\bullet\to N^\vee\]
is given by product with $\beta(E,s)$. \qed
\end{Proof}

\section{Applications to Hodge theory}
\label{dual_deRham}

In this section we prove our main theorems, Theorem~\ref{deRham} and
Theorem~\ref{statements}. We also discuss generalizations of these
results to twisted de Rham complexes.  \bigskip

\noindent
{\bf \large The de Rham complex as a twisted derived intersection}
\bigskip

\noindent
The following is a restatement of Theorem~\ref{deRham}.

\begin{Theorem}
\label{thm:drisdi}
Let $X$ be a smooth scheme over a perfect field $\bbk$ of characteristic
$p>0$.  Let $X'$ be its Frobenius twist, and let $\bS = (T^*X', D)$
denote the Azumaya space whose underlying space is the cotangent
bundle $T^*X'$ of $X'$.  Consider the composite morphism
\[ i': X' = (X',\OO_{X'})\xrightarrow{m}(X',D|_{X'})\xrightarrow{i_D} \bS =  (T^*X',D). \]
Then there are natural isomorphisms in $\D(X')$
\[ F_*\Omega^{\sbt}_X\iso i'^!i'_*\OO_{X'} \iso \left ( i'^*i'_*
  \cO_{X'} \right )^\chk\!\!\!.\] 
Thus $F_*\Omega^\sbt_X$ is formal if and only if $i'^*i'_*\cO_{X'}$
is.
\end{Theorem}

\paragraph
\label{rem:slisdr}
{\bf Remark.} In geometric language we can reformulate the above
theorem by saying that the line bundle associated to the twisted
derived self-intersection problem of $X'$ inside $\bS$ is $\sL =
(F_*\Omega^{\sbt}_X)^{\chk}$.  The complex $F_*\Omega^\sbt_X$ is
formal if and only if $\sL$ is trivial.
\medskip

\noindent
We begin with an easy lemma.

\begin{Lemma}
\label{affine}
Let $f: X\rightarrow Y$ be an affine morphism of schemes and let $\sA$
be a quasi-coherent sheaf of algebras on $X$. For any quasi-coherent sheaf of $\sA$-modules $M$,
its direct image $f_*M$ is a quasi-coherent sheaf of
$f_*\sA$-modules. The functor $f_*$ gives an equivalence
between the category of quasi-coherent $\sA$-modules on $X$ and the category of quasi-coherent $f_*\cA$-modules
on $Y$.
It also gives an equivalence between the corresponding derived
categories. \qed
\end{Lemma}
\medskip

\noindent
{\em Proof of Theorem~\ref{deRham}.} 
Recall from~\eqref{D and DX} that the sheaf of Azumaya
algebras $D$ on $T^*X'$ satisfies $\phi_*D=F_*D_X$, where
$\phi:T^*X' \ra X'$ is the bundle map. We now have 
\begin{align*}
F_*\Omega_X^{\sbt} & = F_*\RsHom_{(X, D_X)}(\cO_X, \cO_X)\\
& =\RsHom_{(X',F_*D_X)}(F_*\OO_X, F_*\OO_X) \\
& =\RsHom_{(X',\phi_*D)}(F_*\OO_X,F_*\OO_X), \\
& = \RsHom_{(X',\phi_*D)}(\phi_*i_{D,*}F_*\OO_X,\phi_*i_{D,*}F_*\OO_X),\\
& = \phi_{*}\RsHom_{(T^*X',D)}(i_{D,*}F_*\OO_X,i_{D,*}F_*\OO_X). 
\end{align*}
Here the first equality is from~\eqref{deRham as Ext}, the second is
Lemma~\ref{affine}, the third is the relationship $\phi_*D=F_*D_X$, the
fourth is the identity $\phi \circ i_D = \id$, and the last is
Lemma~\ref{affine} again.

Recall that the functor 
\[ i_{D,*}:\D(X', D|_{X'}) \ra \D(T^*X', D) \] 
admits a right adjoint $i_D^!$, thus we have
\[ \RsHom_{(T^*X', D)}(i_{D,*} F_* \cO_X, i_{D,*} F_* \cO_X) = i_{D,*}
\RsHom_{(X', D|_{X'})}(F_* \cO_X, i_D^! i_{D,*} F_* \cO_X). \]
Because $\phi\circ i_D=\id$ we conclude that
$$F_*\Omega_X^{\sbt}=\RsHom_{(X',D|_{X'})}
(F_*\OO_X,i_D^!i_{D,*}F_*\OO_X).$$

At this point recall that the sheaf of Azumaya algebras $D|_{X'}$
splits with splitting module $F_*\OO_X$. In fact $F_* \cO_X = m_*
\cO_{X'}$ under the Morita equivalence $\cO_{X'}$-$\gMod
\stackrel{m_*}{\lra} D|_{X'}$-$\gMod$ of~(\ref{Morita}), and the
inverse equivalence is given by the functor 
\[ m^*(\,-\,) = \RsHom_{(X', D|_{X'})}(F_* \cO_X, \,-\,). \]
Moreover, $m^*$ is both a left and a right adjoint to $m_*$ because
the latter functor is an equivalence.  This justifies us to write
either $m^*$ or $m^!$ for the functor above.
 
With these considerations in mind the result of the above calculation
can be rewritten as
\begin{align*}
F_* \Omega_X^\sbt & =
\RsHom_{(X',D|_{X'})}(F_*\OO_X,i_D^!i_{D,*}F_*\OO_X)=m^!i_D^!i_{D,*}m_*\OO_{X'} \\
& = i'^! i'_* \cO_{X'}. 
\end{align*}
The statement that $i'^! i'_* \cO_{X'} \iso \left ( i'^* i'_* \cO_{X'} \right
)^\chk$ is a standard form of duality. Finally, the discussion
in~(\ref{subsec:lineb}) implies that the associated line bundle of
this twisted derived intersection is
$(F_*\Omega^{\sbt}_X)^{\chk}$.\qed 

\paragraph
{\bf Remark.}  The algebra $i^!i_* \cO_{X'}$ is always formal, i.e,
\[ i^!i_* \cO_{X'} \iso \S_{X'}(\Omega_{X'}^1[-1]). \] Since
$i'^!i'_*\cO_{X'}$ and $i^!i_*\cO_{X'}$ have the same cohomology
sheaves anyway, formality of $F_*\Omega_X^\sbt$ is equivalent to the
existence of an isomorphism
\[ i'^!i'_* \cO_{X'} \iso i^! i_* \cO_{X'}, \]
or, dually, the existence of an isomorphism
\[ i'^* i'_* \cO_{X'} \iso i^* i_* \cO_{X'}. \]
This justifies our statement that the above formality statements are
equivalent to the triviality of $\sL$.
\bigskip

\noindent
{\bf \large Splitting on the first infinitesimal neighborhood}
\bigskip

\noindent
In this section we investigate under what circumstances the complex
$F_*\Omega^{\sbt}_X = i'^!i'_* \cO_{X'}$ is formal.  We prove
the following theorem, which is essentially Theorem~\ref{statements}.

\begin{Theorem}\label{D:eq:W}
Let $X$ be a smooth variety over a perfect field $\bbk$ of
characteristic $p>\dim X$. Then the following statements are equivalent.
\begin{itemize}
\item[(1)] $X$ lifts to $W_2(\bbk)$.
\item[(2)] We have $\alpha = 0$, where $\alpha\in H^2(X', T_{X'})$ is
  the class of the extension 
$$0\ra \OO_{X'}\ra F_*\OO_X\xrightarrow{F_*d} F_*Z^1\ra \Omega^1_{X'}\ra 0.$$
Here $Z^1$ denotes the image of $d$ in $\Omega^1_X$.
\item[(3)] The map $i'$ splits to first order.
\item[(4)] $D$ splits on the first infinitesimal neighborhood of $X'$ in
$T^*X'$.
\item[(5)] The line bundle $\sL$ associated to the twisted derived
  self-intersection problem of $X'$ inside $\bS$ is trivial. 
\item[(6)] $F_*\Omega^{\sbt}_X= i'^!i'_*\OO_{X'}$ is formal in $\D(X')$,
  that is, there exists an isomorphism
\[ F_*\Omega^{\sbt}_X= i'^!i'_*\OO_{X'}\iso
\S(\Omega^1_{X'}[-1])=i^!i_*\OO_{X'}. \]
\end{itemize}
\end{Theorem}

\begin{Proof}
We apply Theorem~\ref{thm:twiderint} to the problem of studying the
formality of the twisted derived self-intersection of $X'$ inside $\bS
= (T^*X', D)$.  Conditions (a) and (c) of Theorem~\ref{thm:twiderint}
are automatically satisfied.  Condition (b) is $(3)$ above, and the
theorem thus shows that $(3)$ is equivalent to $(6)$.
\medskip

\noindent
Remark~\ref{rem:anysplit} shows that $(3)$ is also equivalent to
$(4)$. 
\medskip

\noindent
The equivalence of $(5)$ and $(6)$ follows from Remark~\ref{rem:slisdr}.
\medskip

\noindent
Note that $F_*\Omega_X^\sbt \iso i'^!i'_* \cO_{X'}$, so we have 
\[ \tau^{\leq 1} (F_*\Omega_X^\sbt) \iso \tau^{\leq
  1}(i'^!i'_*\cO_{X'}). \] 
The class of the extension in $(2)$ is the obstruction to splitting
the former; the twisted HKR class $\alpha_{\cO_{X'}}$ is by definition
the obstruction to splitting the latter.  We
conclude that $(2)$ is equivalent to $\alpha_{\cO_{X'}} = 0$.
Moreover, Theorem~\ref{twisted} shows that $\alpha_{\cO_{X'}} = 0$ is
equivalent to $(4)$.  Thus we have proved that $(2)$ is equivalent to $(4)$.
\medskip

\noindent
Finally, the equivalence of conditions $(1)$ and $(2)$ is proved in
the original paper of Deligne-Illusie~\cite{DelIll}.\qed
\end{Proof}
\newpage

\noindent
{\bf \large The twisted de Rham complex}
\bigskip

\noindent
We have seen in the previous section that the (dual of the) Frobenius push-forward
of the de Rham complex can be directly regarded as a line bundle over the derived
self-intersection of the zero section in the twisted cotangent bundle
$(T^*X', D)$.  In this section we argue that the twisted de Rham
complex of Barannikov-Kontsevich can also be regarded in a similar way
as a line bundle over a derived intersection in $(T^*X', D)$, as in
Theorem~\ref{twisteddeRham}.  This allows us to give a geometric
realization of the Hodge-de Rham spectral sequence for categories of
singularities.

\paragraph
We start with a scheme $X$ which is smooth over a perfect field $\bbk$ of characteristic $p>0$. We do not assume anymore that $X$ is proper over $k$. We choose  a regular function $f:X\ra \A^1$ on $X$, sometimes called the
superpotential in physics literature.

From these data we construct the context for a twisted derived
intersection as in Section~\ref{sec:twiderint}. As the ambient space $S$ we
take the total space of the cotangent bundle $T^*X'$ of the
Frobenius twist $X'$. The Azumaya algebra on $S$ is $\cA = D$, the ring
of differential operators on $X$.  As one of the two subvarieties of
$S$ we take the zero section $X'\inj T^*X'$.

Let $f' = \pi^*f$ where $\pi:X'\ra X$ is as in (\ref{frobenius}). 
 We take as the second subvariety of $T^*X'$
the graph $X_f'$ of the exact 1-form $df'$. As an abstract variety, 
$X_f'$ is isomorphic to $X'$, with a different embedding into
$T^*X'$. As above, we denote the embeddings $X'\inj T^*X'$ and 
$X_f'\inj T^*X'$ by $i$ and $j$, respectively.
Notice that the underived intersection (which was
previously denoted by $W^0$) is the critical locus $\Crit f'$ of $f'$.

Locally on $X'$, any vector field $\del$ gives rise to a function on
$T^*X'$, linear in the fibers of the map $T^*X' \ra X'$.
The zero section of $T^*X'$ is locally cut out by the
equations $\del' = 0$ for $\del\in T_{X'}$.  Similarly, the graph of $df'$ is cut
out by the equations $\del'-\del'(f') = 0$.

\paragraph
Denote by $\psi$ the automorphism of the ring of differential operators $D_X$ on $X$ that acts on
functions $\cO_X\subset D_X$ as identity and on vector fields as
$\del \mapsto \del - \del f$.  Its
action on the center $\cO_{T^*X'}$ of $D$ is made explicit by
the following lemma.

\begin{Lemma}
\label{lem:actionpsi}
  The automorphism $s$ of $\cO_{T^*X'}$ induced by
  $\psi$ acts on functions as identity and on vector fields $\del'\in T_{X'}$ as
  $\del'\mapsto \del'-\del'(f')$.
\end{Lemma} 

\begin{Proof} 
As discussed in~(\ref{center_D}), $\pi$ induces an isomorphism between
 the sheaves of vector fields on $X$ and on $X'$. For a local vector field $\del'$ on $X'$, we write
$\del$ for the corresponding local vector field on $X$.  A straightforward
calculation shows that 
\begin{align*}
\psi(g^p)&=g^p\\
\psi(\del^p)&=\del^p-\del^{[p]}f-(\del f)^p
\end{align*} 
for any $g\in\cO_X$ and $\del\in T_X$. 
Therefore,
\[\psi(\del^p-\del^{[p]})=\del^p-\del^{[p]}f-(\del f)^p-(\del^{[p]}-\del^{[p]}f)=\del^p-\del^{[p]}-(\del f)^p,\]
where $\del^p-\del^{[p]}$ and $\del^p-\del^{[p]}-(\del f)^p$ are
exactly the images of $\del'$ and $\del'-\del'(f')$ under the
identification $\bar{\iota}:\cO_{T^*X'}\ra Z(D_X)$. \qed
\end{Proof}

\paragraph
\label{def:L}
We shall write $L$ for the $D_X$-module $\psi_*\cO_X$.  As an
$\cO_X$-module it is isomorphic to $\cO_X$, however the action of
$\del \in T_X$ is modified so that $\del$ acts on a local section $g$
of $\cO_X$ by
\[ \del.g = \del(g) - \del f. \] 
We shall think of $\psi$ as an automorphism of the twisted space
$(T^*X', D)$, acting by $s$ on the underlying space $T^*X'$ and by
$\psi$ on $D$.  By Lemma~\ref{lem:actionpsi}, we have $s(X')=X_f'$. 
Since $F_*\cO_X$ splits the restriction $i^*D$ of $D$, $F_*L = F_* \psi_*\cO_X$ splits $j^*D$. 

\paragraph 
We conclude that we are in the context of Section~\ref{sec:twiderint}:
we have an ambient space $\bS = (T^*X', D)$, with two subvarieties
$X', X'_f$ such that $D$ splits on both of them.  We will use the
notations of Section~\ref{sec:twiderint} from now on.  In particular
we shall denote by $i'$, $j'$ the twisted embeddings of $X'$, $X'_f$
into $\bS$, and by $W$ and $W'$ the untwisted and twisted derived
intersections of $X'$ and $X'_f$ (inside $S$ and $\bS$).

In the remainder of this section, we show
that the structure complex of $W$ is the dual of $\Omega^\sbt_{X', \wedge df'}$,
and that the associated line bundle $\sL$ for the above intersection
problem is the dual of $F_*\Omega^\sbt_{X, d-\wedge df}$.

\paragraph 
To explain the above relationship we need to remind the reader about
the two complexes that appear in Theorem~\ref{twisteddeRham}.  The
first is the de Rham complex of the $D_X$-module $L$ given by the
complex of sheaves on $X$
\begin{align*}
\dR(L) = \Omega^\sbt_{X, d-\wedge df}  &= 0 \lra \Omega^0_X
\xrightarrow{d-\wedge df} \Omega^1_X \xrightarrow{d-\wedge df}
\Omega^2_X \lra \cdots;
\intertext{we view it as a 
twisted analogue of the
de Rham complex of $X$. The second is the analogue of the formal complex
  $\S(\Omega^1_{X'}[-1])$; it is the complex of coherent sheaves on $X'$ given by}
\Omega^\sbt_{X', \wedge df'} & = 0 \lra \Omega^0_{X'}
\xrightarrow{\wedge df'} \Omega^1_{X'} \xrightarrow{\wedge df'}
\Omega^2_{X'} \lra \cdots.
\end{align*}

\noindent
The following theorem rephrases Theorem~\ref{twisteddeRham}. 

\begin{Theorem}
\label{thm:twiderha}
Let $W$ and $W'$ be the derived intersections of the zero
section and the graph of $df'$ inside $(T^*X',
\cO)$ and $(T^*X', D)$, respectively.  Then we
have
\begin{itemize}
\item[(a)] the structure complex $\cO_W$ of $W$ is the dual (over $X'$
  ) of $\Omega^\sbt_{X', \wedge df'}$; and
\item[(b)] the associated line bundle $\sL$ of the twisted derived
  intersection problem is the dual (again, over $X'$) of
  $F_*\Omega^\sbt_{X, d-\wedge df}$.
\end{itemize}
If furthermore $\sL$ is trivial, then the two complexes
$\Omega^\sbt_{X', \wedge df'} $ and $F_*\Omega^\sbt_{X, d-\wedge df}$
are quasi-isomorphic.
\end{Theorem}

\begin{Proof}
Part (a) is an easy calculation with an explicit Koszul resolution of
$j_* \cO_{X'}$, so we will concentrate on part (b).  We have the
following sequence of equalities
\begin{align*}F_*\Omega^{\sbt}_{X,d-\wedge df}&=F_*\RsHom^{\sbt}_{(X,D_X)}(\OO_X,L)=\\
  &=\RsHom^{\sbt}_{(X',F_*D_X)}(F_*\OO_X,F_*L)=\\
  &=\RsHom^{\sbt}_{(X',\phi_*D)}(F_*\OO_X,F_*L),
\end{align*} 
where the first equality follows from \eqref{deRham as Ext}, the
second is Lemma \ref{affine} and the third is $F_*D_X=\phi_*D$.

Following the notation of Section~\ref{sec:Azumaya} we denote by $i_D$, $j_D$ the
embeddings of the twisted spaces $(X', i^*D)$, $(X', j^*D)$
into $(T^*X', D)$.  Note that $\phi\circ j_D=\id = \phi\circ
i_D$ because both $i$ and $j$ are sections of the bundle map
$\phi:T^*X'\ra X'$.  We have
\begin{align*}
F_*\Omega^{\sbt}_{X,d-\wedge df}&=\RsHom^{\sbt}_{(X',\phi_*D)}(F_*\OO_X,F_*L)=\\
&=\RsHom^{\sbt}_{(X',\phi_*D)}(\phi_*i_{D,*}F_*\OO_X,\phi_*j_{D,*}F_*L)=\\
&=\phi_*\RsHom^{\sbt}_{(T^*X',D)}(i_{D,*}F_*\OO_X,j_{D,*}F_*L)=\\
&=\phi_*i_*\RsHom^{\sbt}_{(X',i^*D)}(F_*\OO_X,i_D^!j_{D,*}F_*L)=\\
&=\RsHom^{\sbt}_{(X',i^*D)}(F_*\OO_X,i_D^!j_{D,*}F_*L),
\end{align*}
where the second equality comes from $\phi\circ j=\id = \phi\circ i$,
the third is Lemma \ref{affine}, the fourth is adjunction, and the fifth
is again $\phi\circ i=\id$. Since $F_*\cO_X$ and $F_*L$ are splitting
modules of $i^*D$ and $j^*D$ respectively, a similar calculation to
the calculation done in the proof of Theorem~\ref{thm:drisdi} shows
that
\[F_*\Omega^{\sbt}_{X,d-\wedge df} = m^!i_D^!j_{D,*}n_*\OO_{X'}=i'^!j'_*\OO_{X'}.\]
This implies that the dual of $F_*\Omega^{\sbt}_{X,d-\wedge df}$ is
the associated line bundle of the twisted derived intersection. (We
have denoted by $n$ the isomorphism $(X', \cO_{X'}) \iso (X', j^*D)$
induced by the splitting module $F_*L$ of $j^*D$.)\qed
\end{Proof}

\paragraph
Recall that in the case of the untwisted de Rham complex the main
result was a complete characterization of when the associated line
bundle of the intersection is trivial. Specifically, we argued that
$F_*\Omega_X^\sbt$ and $\S(\Omega_{X'}^1[-1])$ are isomorphic if and
only if $X$ lifts to $W_2(\bbk)$.

In the case of the twisted de Rham complex, the natural question is 
whether there exists a quasi-isomorphism
\[
F_*\Omega_{X, d-\wedge df}^\sbt \iso \Omega^{\sbt}_{X', \wedge df'}.\]
Unfortunately, it seems hard to give a precise criterion for the existence of
such an isomorphism. 

Note however that Ogus and Vologodsky prove (see Proof of Theorem~4.23
in \cite{OguVol}) that if we start with a variety $X$ and a proper map
$f:X\to\A^1$ defined over a field of characteristic zero, and consider
reductions of this pair $(X,f)$ to characteristic $p$, then the above
quasi-isomorphism exists for all $p$ that are large enough. They then
use this result to derive the Barannikov-Kontsevich Theorem.

\paragraph 
In this paper, we consider a slightly different question: we ask for
conditions when both complexes $F_*\Omega_{X, d-\wedge df}^\sbt$ and
$\Omega^{\sbt}_{X', \wedge df'}$ are formal.  This turns out to be a
special case of the above problem: these two complexes have isomorphic
cohomology sheaves anyway; therefore, once they are known to be
formal, they are automatically quasi-isomorphic.

\paragraph
As a particular case, the untwisted de Rham theorem corresponds to
$f=0$. Indeed, if $f=0$ the complex
\[\Omega^{\sbt}_{X', \wedge df'}\simeq\S(\Omega_{X'}^1[-1])\]
is clearly formal, and its cohomology sheaves are isomorphic to the
cohomology sheaves as $F_*\Omega_X^\sbt$. Therefore, if $f=0$, the two
complexes are quasi-isomorphic if and only if they are both formal.
\medskip

\noindent
The theorem below is our main result concerning the formality of
twisted de Rham complexes.

\begin{Theorem}
\label{thm:twisteddR1}
Let $X$ be a smooth variety over a perfect field of characteristic $p$, and let
$f:X\to\A^1$ be a morphism.  Assume that the following conditions hold:
\begin{itemize}
  \item[1.] $p>\dim(X)$.
  \item[2.] Both $X$ and $f$ lift to $W_2(\bbk)$.
  \item[3.] The critical locus $Z = \Crit f$ is scheme-theoretically smooth.
  \item[4.] $f$ is constant on each connected component of $Z$.
  \item[5.] The embedding $Z \hookrightarrow X$ is split to first
    order; in other words the short exact sequence below splits
    \[ 0 \ra T_Z\ra T_X|_Z \ra N_{Z/X} \ra 0. \]
\end{itemize}
Then the complexes $F_*\Omega_{X, d-\wedge df}^\sbt$ and
$\Omega^{\sbt}_{X', \wedge df'}$ are both formal and quasi-isomorphic.
\end{Theorem}

\paragraph
{\bf Remark.} If conditions 1 and 3 in the statement of the theorem
are satisfied, formality of $\Omega^{\sbt}_{X', \wedge df'}$ implies
condition 5 by~\cite[Theorem 1.8]{AriCalHab}, thus providing a partial
converse of Theorem~\ref{thm:twisteddR1}.

\begin{Corollary} \label{co:twisteddR}
Suppose that $X$ and $f$ satisfy the conditions 1--5 of Theorem~\ref{thm:twisteddR1}.
Then there
are isomorphisms
\[ \R^i\Gamma(X, \Omega_{X, d-\wedge df}^\sbt) \iso \R^i\Gamma(X,
\Omega^{\sbt}_{X', \wedge df'})\iso \bigoplus_{i-c=p+q}
H^p(Z',\Omega^q_{Z'}\otimes \omega)\] 
where $c$ denotes the codimension of $Z$ in $X$ and $\omega$ denotes
the relative dualizing complex of the embedding $Z'\inj X'$.
\end{Corollary}

\begin{Proof}
Indeed, it is easy to see that since $Z'$ is smooth, the cohomology
sheaves of the complex $\Omega^{\sbt}_{X', \wedge df'}$ are given by
\[H^i(\Omega^{\sbt}_{X', \wedge df'})=\Omega^{i-c}_{Z'}\otimes\omega.\tag*{\qed}\]
\end{Proof}
\medskip

\noindent
Corollary~\ref{co:twisteddR} can now be transferred to characteristic
zero, giving Theorem~\ref{thm:twisteddR}.  For the reader's
convenience we repeat the statement of Theorem~\ref{thm:twisteddR}
below.

\begin{Theorem}
Let $X$ be a smooth variety over a field $\bbk$ of characteristic
zero, and let $f:X\to\A^1$ be a proper morphism. Suppose that the
following conditions are satisfied:
\begin{itemize}
  \item[1.] The critical locus $Z = \Crit f$ is scheme-theoretically smooth.
  \item[2.] The embedding $Z \hookrightarrow X$ is split to first order, that
    is, the short exact sequence below splits
    \[ 0 \ra T_Z\ra T_X|_Z \ra N_{Z/X} \ra 0. \]
\end{itemize}

Then there exist isomorphisms
\[ \R^i\Gamma(X, \Omega_{X, d-\wedge df}^\sbt) \iso \R^i\Gamma(X,
\Omega^{\sbt}_{X, \wedge df})\iso \bigoplus_{i-c=p+q}
H^p(Z,\Omega^q_{Z}\otimes \omega)\] 
where $c$ denotes the codimension of $Z$ in $X$ and $\omega$ denotes
the relative dualizing complex of the embedding $Z\inj X$.
\end{Theorem}

\begin{Proof} 
The theorem is deduced from the characteristic $p$ statement
(Corollary~\ref{co:twisteddR}) in a manner parallel to the argument of
\cite{OguVol}.  Note that this argument relies on a very subtle
statement concerning specialization of the de Rham cohomology to
characteristic $p$, see the first claim of
\cite[Theorem~4.23]{OguVol}.
\qed
\end{Proof}

\paragraph {\bf Remark.} The isomorphism
\[\R^i\Gamma(X,
\Omega^{\sbt}_{X, \wedge df})\iso \bigoplus_{i-c=p+q}
H^p(Z,\Omega^q_{Z}\otimes \omega)\]
follows directly from \cite[Theorem~1.8]{AriCalHab} (which is the untwisted version of Theorem~\ref{thm:twiderint}).
Thus, Theorem~\ref{thm:twisteddR} is also implied by the combination of \cite[Theorem~1.8]{AriCalHab}
and the Barannikov-Kontsevich Theorem (Theorem~\ref{BarKon}).

\paragraph 
Let us now prove Theorem~\ref{thm:twisteddR1}. Its proof is a
straightforward application of Theorem~\ref{thm:twiderint} (which in
turn generalizes the main result of~\cite{AriCalHab}, similar in
nature also to the main result of~\cite{Gri}).  The only non-trivial
aspect of the proof is checking condition (a) of
Theorem~\ref{thm:twiderint}; that is, we need to verify that the
splitting modules of the Azumaya algebra $D$ along $X'$ and $X'_f$
agree along the (non-derived) intersection $Z' = X' \cap X_f'$. This
requires a study of exponential functions in positive characteristic.

\paragraph
Let $J\subset\cO_X$ be the Jacobian ideal of $f$, that is, the ideal of the
closed subscheme $Z\subset X$. Put 
\[J^{[p]}=\langle g^p|g\in J \rangle,\]
so that $J^{[p]}$ is the ideal of the closed subscheme $Z^{[p]} = F^{-1}(Z')\subset X$.
Note that $J^{[p]}\subset\cO_X$ is a sub-$D_X$-module.

Recall that $D$ is trivialized along the zero section $X'$ of $T^*X'$
by the module $E = F_* \cO_X$. By Lemma~\ref{lem:actionpsi}, acting
by $\psi$ on $D$ moves the zero section $X'$ to the section $X'_f$,
and hence the module $E_f = F_* \psi_* \cO_X = F_* L$ is a splitting
module of $D$ along $X_f'$.  To complete the proof of Theorem~\ref{thm:twisteddR1},
we need to show that $E|_{Z'} \iso
E_f|_{Z'}$ {\em as $F_*D|_{Z'}$-modules}.  Since the map $F:X \ra X'$ is
affine, this is equivalent to finding an isomorphism of $D_X$-modules
\[\cO_X/J^{[p]}\simeq L/L\cdot J^{[p]}.\]

Recall that $L$ is a trivial rank one $\cO_X$-module, but derivations $\del\in T_X\subset D_X$ act
on it by $\del-\del f$. Thus, an isomorphism 
\[\cO_X/J^{[p]}\simeq L/L\cdot J^{[p]}\]
is given by multiplication by an invertible function
\[g\in(\cO_X/J^{[p]})^\times\]
that must satisfy the differential equation
\[ d\log(g) = df \bmod J^{[p]}; \]
as usual, $d\log(g)=dg/g$. 

We will call a solution $g$ to the above differential equation an {\em
  approximate exponential} of $f$. Observe that in characteristic $p$ the usual power series for
$\exp(f)$ is not defined past the $p$-th term.

\paragraph
At this point, the existence of an approximate
exponential follows from the work of Ogus and Vologodsky~\cite[Proposition
4.28]{OguVol}. This completes the proof of
Theorem~\ref{thm:twisteddR1}. Since the relevant part of \cite{OguVol} is spread over several
claims, we restate their result and sketch its proof below for the reader's convenience.

\begin{Proposition}\label{pp:exponential}
As before, let $X$ be a smooth variety over a perfect field $\bbk$ of characteristic $p$, and let
$f:X\to\A^1$ be a morphism. Denote by $Z=\Crit f$ the critical locus of $f$, and let $J\subset\cO_X$
be the Jacobian ideal of $f$.  
Suppose that the following conditions hold:
\begin{itemize}
  \item[1.] Both $X$ and $f$ lift to $W_2(\bbk)$.
  \item[2.] $f$ is constant on each reduced connected component of $Z$.
  \item[3.] $p>cn$, where $c$ is the maximal codimension of $Z$ in $X$, and $n$ is the smallest number such that $(\sqrt{J})^n\subset J$.
  
\end{itemize}
Under these assumptions, an approximate exponential exists: there is $g\in(\cO_X/J^{[p]})^\times$ such that
\[ d\log(g) = df \bmod J^{[p]}. \]
\end{Proposition}

\begin{Proof} Without loss of generality, we may assume that $Z$ is connected. Subtracting a constant from $f$
we may assume that $f\in\sqrt{J}$. 

Define $\alpha\in\cO_X/J^{[p]}$ as follows. Locally on $X$ there exists a lift $\tilde F$ of the Frobenius map $F$ to
$W_2(\bbk)$; put 
\[ \alpha = \frac{1}{p}\left ( \tF(\tf) - \tf^p \right) \bmod p, \]
where $\tf$ is the lift of $f$. One can check that $\alpha$ does not depend on the choice of $\tilde F$, and therefore
the local expressions glue to define $\alpha$ globally.

Now define $g$ by the formula
\[ g = \AH(f) \cdot \left ( 1 + \alpha +\frac{\alpha^2}{2!} + \cdots +
  \frac{\alpha^{p-1}}{(p-1)!} \right ). \]
Here 
\[ \AH(x) = \exp\left (x +\frac{x^p}{p} + \frac{x^{p^2}}{p^2} + \cdots
\right ) \]  
  is the Artin-Hasse exponential~\cite{ArtHas}.
It is a formal power series whose coefficients are rational numbers without $p$'s
in the denominator, and therefore $\AH(x)$ makes sense over fields of characteristic
$p$.

It is not hard to check that $g$ is indeed an approximate exponential. Indeed, 
\[d\log(\AH(f))=(1+f^{p-1}+f^{p^2-1}+\dots)df =(1+f^{p-1})df\bmod J^{[p]},\]
because $f^{p^2-1}\in J^{[p]}$. Local calculation shows that $\alpha^c\in(\sqrt{J})^{[p]}$, therefore $\alpha^{p-1}=0\bmod J^{[p]}$,
which implies that
\[d\log\left ( 1 + \alpha +\frac{\alpha^2}{2!} + \cdots +
  \frac{\alpha^{p-1}}{(p-1)!} \right )=d\alpha\bmod J^{[p]}.\]
It remains to notice that $d\alpha=f^{p-1}df$, and hence 
\begin{align*}
d\log(g) &=d\log(\AH(f))+d\log\left ( 1 + \alpha +\frac{\alpha^2}{2!} + \cdots +
  \frac{\alpha^{p-1}}{(p-1)!} \right ) \\
 & =df\bmod J^{[p]}
\end{align*}
  as claimed.\qed
\end{Proof}

\paragraph {\bf Remark.} 
\label{subsec:defh} 
We conclude the paper with an observation concerning approximate
exponentials. As we have already mentioned, $\alpha^{p-1}=0\bmod J^{[p]}$
under the hypotheses of Proposition~\ref{pp:exponential}. For this
reason, we can rewrite the formula for $g$ as
\[g=\AH(f)\cdot\AH(\alpha).\]
This suggests a different informal approach to finding an approximate
exponential as an infinite product. Namely, suppose that there exist
lifts of $\tX$ and $\tf$ not just to the ring of second Witt vectors
$W_2(\bbk)$, but to the ring of infinite-level Witt vectors
$W(\bbk)$. We can now consider lifts $\tF$ of the Frobenius map $F$ to
$W(\bbk)$; generally speaking, such $\tF$ exists only locally.

Consider the sequence
\begin{align*}
h_0 & = \tf \\
h_1 & =  \frac{1}{p}\left ( \tF(\tf) - h_0^p \right) \\
h_2 & =  \frac{1}{p^2} \left ( \tF^2(\tf) - h_0^{p^2} -p\cdot h_1^p \right
) \\
& \vdots  \\
h_n & =  \frac{1}{p^n} \left ( \tF^n(\tf) -\sum_{j=0}^{n-1} p^j
  h_j^{p^{n-j}} \right ) \\
& \vdots 
\end{align*}
These expressions are well defined over $W(\bbk)$ (in other words the
expressions in parentheses are divisible by the corresponding power of
$p$). One may consider as a candidate for an
approximate exponential of $f$ the reduction modulo $p$ of the infinite product
\[
\tilde g=\AH(h_0) \cdot \AH(h_1) \cdot \AH(h_2)\cdots.
\]
(Here we are completely ignoring convergence issues.) 
The approximate exponential used in Proposition~\ref{pp:exponential} 
is the reduction of the truncated product $\AH(h_0)\cdot \AH(h_1)$.

The recursive expressions for $h_i$ imply the identity
\[\tilde g=\exp\left ( \tf + \frac{\tF(\tf)}{p} + \frac{\tF^2(\tf)}{p^2}
  + \cdots \right )\]
  over $W(\bbk)[p^{-1}]$. Taking the derivative of both sides, we see that
  \[d\tilde g=d\tf+\frac{d\tF(\tf)}{p}+\frac{d\tF^2(\tf)}{p^2}+\cdots.\]
  It is easy to see that $d\tF^k(\tf)$ is actually divisible by $p^k$ over $W(\bbk)$, 
  so the right-hand side of the formula makes sense over $W(\bbk)$. Moreover, 
  the reduction of $p^{-k}d\tF^k(\tf)$ modulo $p$ vanishes modulo $J^{[p]}$. This justifies
  the expectation that $\tilde g$ is an approximate exponential of $\tilde f$.

The above formulas are perhaps more conceptual than the truncated expression used in Proposition~\ref{pp:exponential}. Unfortunately, 
there are two issues with this approach, and it is not clear how to make it 
rigorous. First of all, the above infinite products and sums may diverge.
Also, even when the expression for $\tilde g$ is well defined, it may depend
on the choice of local lift $\tF$.  
In some cases, these issues can be resolved by restricting the class of lifts $\tF$
that one considers.

\end{document}